\documentclass[11pt,reqno]{amsart}
\usepackage{graphicx}
\usepackage{amssymb,latexsym,epsfig,mathrsfs,graphpap,graphics,amsmath,amssymb}
\usepackage{amstext}
\usepackage{verbatim}
\usepackage{color}

\newtheorem{conj.}[thm]{Conjecture}

\theoremstyle{definition}

\theoremstyle{remark}

\numberwithin{equation}{section}

\begin{document}

\begin{flushleft}
  {\bf\Large {Convolution Based Special Affine Wavelet Transform and Associated Multi-resolution Analysis}}
\end{flushleft}

\parindent=0mm \vspace{.4in}

  {\bf{Firdous A. Shah$^{\star}$ and  Waseem Z. Lone$^{\star}$  }}

 \parindent=0mm \vspace{.1in}
{\small \it $^{\star}$Department of  Mathematics,  University of Kashmir, South Campus,
 Anantnag-192101, Jammu and Kashmir, India. E-mail: $\text{fashah@uok.edu.in}$;\,$\text{lwaseem.scholar@kashmiruniversity.net}$}

\parindent=0mm \vspace{.2in}
{\small {\bf Abstract.} In this paper, we study the convolution structure in the special affine Fourier domain to combine the advantages of the well known special affine Fourier and wavelet transforms into a novel integral transform coined as special affine wavelet transform and investigate the associated constant Q-property in the joint time-frequency domain. The preliminary analysis encompasses the derivation of the fundamental properties, orthogonality relation, inversion formula and range theorem. Finally, we extend the scope of the present study by introducing the notion of  multi-resolution analysis associated with special affine wavelet transform and the construction of orthogonal special affine wavelets. We call it special affine multi-resolution analysis. The necessary and sufficient conditions pertaining to special affine Fourier domain under which the integer shifts of a chirp modulated functions form a Riesz basis or an orthonormal basis for a multi-resolution subspace is established.

\parindent=0mm \vspace{.1in}

\parindent=0mm \vspace{.1in}
}

\section{Introduction}
\parindent=0mm \vspace{.0in}

An obligatory addition to the theory of Fourier transforms was presented by Abe and Sheridan \cite{1,2}, in the form of special affine Fourier transform (SAFT). The SAFT is a six parameter class of linear integral transforms which embraces several well known unitary transforms, say the Fourier transform \cite{3,4}, the fractional Fourier transform (FrFT) \cite{5,6}, the linear canonical transform (LCT) \cite{7}, the Fresnel transform \cite{8} and the scaling operations \cite{9}. The SAFT can be regarded as a time-shifting and frequency modulated version of the well known linear canonical transform \cite{10,11,12} and is defined as

\begin{align*}
 \mathrm S_\bold M [f](\omega)=
   \begin{cases}
     &\int_\mathbb R f(t) \, \mathcal K_\bold M (t,\omega) \, dt, \hspace{4.9cm} B \neq 0 \\
     &\sqrt{D} \, \exp \left\{ \frac{i}{2} \big( CD(\omega - p)^2 + 2 \omega q \big) \right\} f \big( d(\omega - p) \big), 		           \quad B=0 ,
   \end{cases}   																							\tag{1.1}
\end{align*}

where $ \bold M $ denotes the $2 \times 3$ unimodular matrix $ \begin{pmatrix}
A & B  & : & p \\
C & D  & : & q
\end{pmatrix} $ and $ \mathcal K_\bold M (t,\omega) $ denotes the kernel of the SAFT given by

\begin{align*}
  \mathcal K_\bold M (t,\omega) = \frac{1}{\sqrt{i 2\pi B}} \, \exp \left\{ \frac{i}{2B} \Big( At^2 + 2t(p-\omega) -    				2\omega(Dp-Bq) + D \big( \omega^2 + P^2 \big) \Big) \right\} .  									\tag{1.2}
\end{align*}

Throughout this article, we shall only consider the case $B \neq 0$, since SAFT is just a chirp multiplication operation in case $B=0$. The inverse SAFT corresponding to $(1.1)$ is given by

\begin{align*}
  f(t) = \int_\mathbb R \mathrm S_\bold M [f](\omega)\,\mathcal K_{\bold M^{-1}} (\omega,t)\, d\omega ,      \tag{1.3}
\end{align*}

where
\begin{align*}
  \bold M^{-1} = \begin{pmatrix}
   					D & -B & : & Bq-Dp \\
  				   -C & A & : & Cp-Aq
				 \end{pmatrix}.
\end{align*}

Moreover, the Parsevals formula for the SAFT is given by

\begin{align*}
  \Big\langle f(t),g(t) \Big\rangle_{L^2(\mathbb R)} = \Big\langle \mathrm S_\bold M [f](\omega), \mathrm S_\bold M [g](\omega) \Big\rangle_{L^2(\mathbb R)} , \quad \forall \, f,g \in L^2(\mathbb R). 									\tag{1.4}
\end{align*}

\parindent=8mm\vspace{.1in}
The birth of the SAFT has attained a respectable status within a short period of time and has been applied to optical, electrical and communication systems, quantum mechanics and several other fields of science and technology \cite{13,SAA,14,15}. Despite the versatility of its applications, SAFT suffers a major drawback due to its global kernel involved in $(1.1$) and therefore is inadequate in situations demanding a joint analysis of time and spectral characteristics of a signal. To circumvent this issue, an immediate concern is to patch up the existing SAFT by adjoining certain excellent localization features. The most appropriate candidates are the wavelet functions, which inherits nice localization properties along with additional characteristic features such as orthogonality, vanishing moments and self adjustability. By intertwining the ideas of SAFT and wavelet transforms, we introduce a novel integral transform coined as, special affine wavelet transform (SAWT) in the context of time-frequency analysis by invoking the convolution structure in the realm of special affine Fourier
transform.

\parindent=8mm\vspace{.1in}
On the other hand multi-resolution analysis (MRA) is an important mathematical tool since it provides a natural framework for understanding and constructing discrete wavelet systems \cite{MM,D}. In recent years there has been a considerable interest in the problem of constructing wavelets via MRA approach namely MRA associated with fractional Fourier transform \cite{SX}, MRA associated with linear canonical transform \cite{JW}. But the concept of multi-resolution analysis associated with the SAWT is yet to be explored. Motivated and inspired by the concept of MRA in $L^2(\mathbb R)$, the notion of MRA in the realm of special affine Fourier domain is introduced and the construction of a Riesz basis or an orthonormal basis is derived.

\parindent=8mm\vspace{.1in}
The rest of the article is structured as follows: In Section 2, we introduce the convolution based special affine wavelet transform and study the associated constant Q-property . In Section 3 , we study the associated fundamental properties. Finally in Section 4, the discrete special affine wavelet transform is proposed, including the definition of special affine multi-resolution analysis and the necessary and sufficient conditions to generate a Riesz basis or an orthonormal basis.

\section{Convolution based special affine wavelet transform and the associated constant Q-property}

\parindent=0mm\vspace{.0in}
In this section, we shall firstly briefly recall the notion of special affine convolution and then formally propose the convolution based special affine wavelet transform. In the sequel we study the constant Q-property associated with the proposed transform.

\parindent=0mm\vspace{.1in}
{\bf Definition 2.1.} For any pair of functions $ f, g \in L^2(\mathbb R)$, the special affine convolution $*_\bold M$ is
defined by \cite{14}

\begin{align*}
  h(t) = \big( f *_\bold M g \big)(t)= \frac{1}{\sqrt{i2\pi B}} \, \exp \left\{ i \frac{D}{2B} p^2 \right\} \int_\mathbb R    	f(\tau) \, g(t-\tau) \, \exp \left\{ i \frac{A\tau}{B} (t-\tau) \right\} d\tau       		      \tag{2.1}
\end{align*}

and the associated convolution theorem reads

\begin{align*}
  \mathrm S_\bold M [h](\omega) = \exp \left\{ \frac{i}{2B} \Big( 2\omega(Dp-Bq) - D\omega^2 \Big) \right\} \mathrm S_   			\bold M [f](\omega) \, \mathrm S_\bold M [g](\omega) .      							                  \tag{2.2}
\end{align*}

Based on the SAFT convolution defined in $(2.1)$, we shall introduce the notion of convolution based special affine wavelet transform.

\parindent=0mm\vspace{.1in}
{\bf Definition 2.2.} For any finite energy signal $f \in L^2(\mathbb R)$, the continuous special affine wavelet transform of $f$ with respect to the wavelet $\psi \in L^2(\mathbb R)$, is defined by

\begin{align*}
  \Big( \mathcal W_{\bf M}^\psi \Big)(b,a) &=  \exp \left\{ -i \frac{A}{2B} \left(1 - \frac{1}{a^2} \right)b^2 \right\} 			\left[ f(t)*_\bold M \frac{1}{\sqrt{a}} \, \overline{\psi \left( \frac{-t}{a} \right)} \exp \left\{ -i \frac{A}{2B} 			t^2 \right\}  \right](b)  \\
  &= \int_\mathbb R f(t) \, \overline{ \psi^{\bold M}_{b,a}(t)} \, dt, \quad b\in \mathbb R, a \in \mathbb R^+ , \tag{2.3}
\end{align*}

where
\begin{align*}
   \psi^{\bold M}_{b,a}(t) =  \frac{1}{\sqrt{i2\pi a B}} \, \psi \left( \frac{t-b}{a} \right) \exp \left\{ - \frac{i}			{2B} \left( At^2 + Dp^2 - A\Big( \frac{b}{a} \Big)^2 \right) \right\}.  							 \tag{2.4}
\end{align*}

\parindent=0mm\vspace{.1in}
As a consequence of Definition 2.2, we have the following deductions

\begin{itemize}
 \item When $\bold M = (A,B,C,D:0,0)$, the special affine wavelet transform $(2.3)$ reduces to the LCT wavelet transform 					\cite{JW}.
 \item For $\bold M = (\cos\theta,\sin\theta,-\sin\theta,\cos\theta:0,0)$, $(2.3)$ reduces to the fractional wavelet transform 				\cite{Fr}.
 \item When $\bold M = (0,1,-1,0:0,0)$, the proposed special affine wavelet transform $(2.3)$ reduces to the classical wavelet  			transform \cite{MM}.
\end{itemize}

\parindent=0mm\vspace{.1in}
{\bf Example 2.3.} Consider the function $f(t) = \exp\left\{-\big( \alpha t + \beta t^2\big) \right\}, \, \alpha ,\beta >0 $ and the Morlet function given below
\begin{align*}
  \psi(t)= \exp\left\{-\left( i\gamma t - \frac{t^2}{2} \right) \right\}.
\end{align*}
	
Then, the translated and scaled versions of  $\psi(t)$ are given by

\begin{align*}
  \psi \left( \frac{t-b}{a} \right)= \exp\left\{ i\gamma \left( \frac{t-b}{a} \right) - \frac{t^2 -2tb +b^2}{2a^2}  \right\}.
\end{align*}

Consequently, the special affine wavelet transform of $f(t)$ and the window function $\psi(t)$ is given by

\begin{align*}
  \mathcal W^{\bf M}_{\psi} \big[ f \big]&= \frac{1}{\sqrt{i2\pi a B}} \int_\mathbb R \exp \Big\{-\big( \alpha t + \beta 				   t^2\big) \Big\} \exp\left\{ -i\gamma \left( \frac{t-b}{a} \right) - \frac{t^2 -2tb +b^2}{2a^2}  \right\}    \\
  & \hspace{4cm} \times\exp \left\{  \frac{i}{2B} \left( At^2 + Dp^2 - A\Big( \frac{b}{a} \Big)^2 \right) \right\} dt    \\
  &= \frac{1}{\sqrt{i2\pi a B}} \exp \left\{  \frac{i}{2B} \left( Dp^2 - A\Big( \frac{b}{a} \Big)^2 \right) - \frac{1}{2} \Big( 		\frac{b}{a} \Big)^2 + i\gamma\frac{b}{a}  \right\} \int_\mathbb R \exp \Big\{-\big( \alpha t + \beta t^2 \big)\Big\} \\
  & \hspace{4cm} \times \exp\left\{ -it\frac{\gamma}{a} - \frac{t^2 -2tb}{2a^2}  \right\} \exp \left\{  \frac{i}{2B} At^2 						\right\} dt  \\
  &= \frac{1}{\sqrt{i2\pi a B}} \exp \left\{  \frac{i}{2B} \left( Dp^2 - A\Big( \frac{b}{a} \Big)^2 \right) - \frac{1}{2} \Big( 		\frac{b}{a} \Big)^2 + i\gamma\frac{b}{a}  \right\} 		\\
  & \hspace{4cm} \times\int_\mathbb R \exp\left\{ -t^2 \left( \beta -i\frac{A}{2B} + \frac{1}{2a^2} \right) -t \left( \alpha +i 					\frac{\gamma}{a} - \frac{b}{a^2}\right)  \right\} dt    \\
  &= \sqrt{\frac{a}{a^2A + i2a^2\beta B +iB}} \, \exp \left\{  \frac{i}{2B} \left( Dp^2 - A\Big( \frac{b}{a} \Big)^2 \right) - 				\frac{1}{2}\Big( \frac{b}{a} \Big)^2 + i\gamma\frac{b}{a}  \right\}    \\
  & \hspace{4cm} \times \exp \left\{ \frac{B\big( b - i\gamma a - \alpha a^2 \big)^2}{2a^2 \big( 2a^2\beta B -i a^2A +B \big)}  			\right\}.
\end{align*}

\parindent=8mm\vspace{.1in}
We now intend to derive a fundamental relationship between the SAFT $(1.1)$ and the proposed special affine wavelet transform $(2.3)$. With the aid of this expression, we shall study the Q-property of the proposed transform and obtain the associated joint time-frequency resolution.

\parindent=0mm\vspace{.1in}
{\bf Proposition 2.4.} {\it Let $\Big( \mathcal W^\psi_{\bf M} \Big)(b,a)$ and $ \mathrm S_{\bf M}[f](\omega)$ be the continuous special affine wavelet transform and the special affine Fourier transform of any finite energy signal $f\in L^2(\mathbb R)$. Then the expression $(2.3)$ in the SAFT domain is given by

\begin{align*}
  \Big( \mathcal W^\psi_{\bf M} f \Big)(b,a) &= \sqrt{a}\, \exp \left\{ -i \frac{A}{2B} \Big( 1 - \frac{1}{a^2} \Big) b^2    	\right\} \int_\mathbb R \exp \left\{ \frac{i}{2B} \Big( 2a\omega(Dp-Bq) - Da^2\omega^2 \Big) \right\}  \\
  &\quad\times \, \mathrm S_{\bf M}[f](\omega) \, \mathrm S_{\bf M} \left[ \exp \left\{ \frac{i}{2B} \Big( 2tp(a-1) - 			At^2 \Big) \right\} \overline{\psi(-t)}\right](a\omega) \, \mathcal K_{\bf M^{-1}}(t,\omega) \, d\omega ,  	\tag{2.5}
\end{align*}

where $ \mathcal K_{\bf M}(t,\omega) $ is given by $(1.2)$.
}

\parindent=0mm\vspace{.1in}
{\it Proof.} Applying the definition of special affine Fourier transform $(1.1)$ together with the implication of $(2.1)$, we have

\begin{align*}
  &\mathrm S_{\bf M} \left[ \Big( \mathcal W^\psi_{\bf M} f \Big)(b,a) \right](\omega)  \\
   &\quad= \exp \left\{ -i \frac{A}{2B} \Big(1 - \frac{1}{a^2} \Big)b^2 \right\} \, \mathrm S_{\bf M} \left[ \left( 					f(t)*_\bold M \frac{1}{\sqrt{a}} \, \overline{\psi \left( \frac{-t}{a} \right)} \exp \left\{ -i \frac{A}{2B}t^2 				\right\} \right)(b) \right](\omega) 	\\
   &\quad= \exp \left\{ -i \frac{A}{2B} \Big(1 - \frac{1}{a^2} \Big)b^2 \right\} \exp \left\{ \frac{i}{2B} \Big( 				2\omega(Dp-Bq) - D\omega^2 \Big) \right\} 		\\
   &\hspace{4cm} \times \mathrm S_{\bf M}[f](\omega) \, \mathrm S_{\bf M} \left[ \frac{1}{\sqrt{a}} \, \overline{\psi \left( 				\frac{-t}{a} \right)} \exp \left\{ -i\frac{A}{2B}t^2 \right\} \right](\omega).                \tag{2.6}
\end{align*}

We have,
\begin{align*}
  &\mathrm S_{\bf M} \left[ \frac{1}{\sqrt{a}} \, \overline{\psi \left( \frac{-t}{a} \right)} \exp \left\{ -i\frac{A}{2B}			t^2 \right\} \right](\omega)	\\
  &\quad= \frac{1}{\sqrt{a}} \int_\mathbb R \overline{\psi \left( \frac{-t}{a} \right)} \exp \left\{ -i\frac{A}{2B} t^2 			\right\} \mathcal K(t,\omega)\, d\omega		\\
  &\quad=  \frac{1}{\sqrt{i2\pi aB}} \int_\mathbb R \overline{\psi \left( \frac{-t}{a} \right)} \exp \left\{ \frac{i}{2B}  	\Big( 2t(p-\omega) - 2\omega(Dp-Bq) + D \big( \omega^2+p^2 \big) \Big) \right\} dt		\\
  &\quad=  \sqrt{\frac{a}{i2\pi B}} \,\, \exp \left\{ \frac{i}{2B} \Big( D \big(\omega^2+p^2 \big) - 2\omega(Dp-Bq) \Big) 		\right\} \int_\mathbb R \overline{\psi(-z)} \, \exp \left\{ \frac{i}{2B} \Big( 2az(p-\omega) \Big) \right\} dz	 \\
  &\quad=  \sqrt{\frac{a}{i2\pi B}} \,\, \exp \left\{ \frac{i}{2B} \Big( D \big(\omega^2+p^2 \big) - 2\omega(Dp-Bq) \Big) 		\right\} \int_\mathbb R \overline{\psi(-z)} \, \exp \left\{ \frac{i}{2B} \Big( 2az(p-\omega) \Big) \right\}  \\
  &\hspace{1.5cm} \times \, \exp \left\{ \frac{i}{2B} \Big( Az^2 + D\big( (a\omega)^2+p^2 \big) + 2z(p-a\omega) - 2a					\omega(Dp-Bq) \Big) \right\}   \\
  &\hspace{1.5cm} \times \, \exp \left\{ -\frac{i}{2B} \Big( Az^2 + D\big( (a\omega)^2+p^2 \big) + 2z(p-a\omega) - 2a					\omega(Dp-Bq) \Big) \right\} dz		
  \end{align*}
\begin{align*}
  &=  \sqrt{a} \, \exp \left\{ \frac{i}{2B} \Big( D \omega^2 \big( 1 - a^2 \big) + 2\omega(Dp-Bq)(a-1) \Big) \right	\}							\hspace{4.5cm}	\\
  &\hspace{1.5cm} \times \int_\mathbb R \exp \left\{ \frac{i}{2B} \Big( 2zp(a-1) - Az^2 \Big) \right\} \overline{\psi(-z)} \, 				\mathcal K_{\bf M} (z,a\omega)\, dz 	\\
  &=  \sqrt{a} \, \exp \left\{ \frac{i}{2B} \Big( D \omega^2 \big( 1 - a^2 \big) + 2\omega(Dp-Bq)(a-1) \Big) \right	\} 		\\
  & \hspace{1.5cm} \times \mathrm S_{\bf M} \left[ \exp \left\{ \frac{i}{2B} \Big( 2zp(a-1) - Az^2 \Big) \right\} 									\overline{\psi(-z)}	\right](a\omega).
\end{align*}

By virtue of above equality, $(2.6)$ can be expressed as

\begin{align*}
&\mathrm S_{\bf M} \left[ \Big( \mathcal W^\psi_{\bf M} f \Big)(b,a) \right](\omega)  \\
   &\quad= \sqrt{a} \, \exp \left\{ -i \frac{A}{2B} \Big(1 - \frac{1}{a^2} \Big)b^2 \right\} \exp \left\{ \frac{i}					{2B} \Big( 2a\omega(Dp-Bq) - Da^2\omega^2 \Big) \right\} 		\\
   &\hspace{3cm} \times  \, \mathrm S_{\bf M}[f](\omega)\, \mathrm S_{\bf M} \left[ \exp \left\{ \frac{i}{2B} \Big( 						 	2zp(a-1) - Az^2 \Big) \right\} \overline{\psi(-z)}	\right](a\omega).
\end{align*}

By implementing inverse SAFT $(1.3)$, the special affine wavelet transform $(2.3)$ can be expressed as

\begin{align*}
  \Big( \mathcal W^\psi_{\bf M} f \Big)(b,a) &= \sqrt{a} \, \exp \left\{ -i \frac{A}{2B} \Big(1 - \frac{1}{a^2} \Big)b^2 				\right\} \int_\mathbb R \exp \left\{ \frac{i}{2B} \Big( 2a\omega(Dp-Bq) - Da^2\omega^2 \Big) \right\} 	\\
  &\quad \times \, \mathrm S_{\bf M}[f](\omega) \, \mathrm S_{\bf M} \left[ \exp \left\{ \frac{i}{2B} \Big( 						 	2zp(a-1) - Az^2 \Big) \right\} \overline{\psi(-z)} \right](a\omega)\, \mathcal K_{\bf M^{-1}} (z,\omega) \, d							\omega.
\end{align*}

This completes the proof of the Proposition $2.4$. \quad\fbox

\parindent=8mm\vspace{.1in}

As a consequence of Proposition $2.4$, we conclude that if the analyzing functions  $\psi_{b,a}^{\bf M} (t)$ are supported in the time-domain or the special affine Fourier domain, then the proposed transform $\mathcal W_{b,a}^\psi (t)$ is accordingly supported in the respective domains. This implies that the special affine wavelet transform is capable of providing the simultaneous information of the time and the special affine frequency in the time-frequency domain. To be more specific,
suppose that  $\psi(t)$ is the window with centre $E_\psi$ and radius $\Delta_\psi$ in the time domain. Then, the centre and radii of the time-domain window function $\psi_{b,a}^{\bf M} (t)$ of the proposed transform $(2.3)$ is given by

\begin{align*}
  E \Big[ \psi_{b,a}^{\bf M} (t) \Big] =  \frac{\int_\mathbb R t \, \big| \psi_{b,a}^{\bf M} (t) \big|^2 dt}{\int_\mathbb R     		\big| \psi_{b,a}^{\bf M} (t) \big|^2 dt} =  \frac{\int_\mathbb R t \, \big| \psi_{b,a} (t) \big|^2 dt}{\int_\mathbb R     				\big| \psi_{b,a} (t) \big|^2 dt} = E \Big[ \psi_{b,a} (t) \Big] = b + a E_\psi           \tag{2.7}
\end{align*}

\parindent=0mm\vspace{.0in}
and
\begin{align*}
  \Delta \Big[ \psi_{b,a}^{\bf M} (t) \Big] &=  \left( \dfrac{\int_\mathbb R \big(t - (b + a E_\psi ) \big) \big| \psi_{b,a}				^{\bf M} (t) \big|^2 dt}{\int_\mathbb R \big| \psi_{b,a}^{\bf M} (t) \big|^2 dt} \right)^{1/2} =  \left( \dfrac{\int_					\mathbb R (t - b - a E_\psi ) \big| \psi_{b,a} (t) \big|^2 dt}{\int_\mathbb R \big| \psi_{b,a} (t) \big|^2 dt} 					 \right)^{1/2}     \\
  &= \Delta \Big[ \psi_{b,a} (t) \Big] = a \Delta_\psi,                                            \tag{2.8}
\end{align*}

respectively. Let $H(\omega)$ be the window function in the special affine Fourier transform domain given by

\begin{align*}
  H(\omega)= \mathrm S_{\bf M} \left[ \exp \left\{ \frac{i}{2B} \Big( 2tp(a-1) - At^2 \Big) \right\} \overline{\psi(-t)}	              					\right](\omega).
\end{align*}

Then, we can derive the center and radius of the special affine Fourier domain window function

\begin{align*}
  H(a\omega)= \mathrm S_{\bf M} \left[ \exp \left\{ \frac{i}{2B} \Big( 2tp(a-1) - At^2 \Big) \right\} \overline{\psi(-t)}							\right](a\omega)
\end{align*}

appearing in $(2.5)$ as

\begin{align*}
  E \Big[ H(a\omega) \Big] = \dfrac{\int_\mathbb R (a\omega)\big|  H(a\omega) \big|^2 d\omega}{\int_\mathbb R \big|  H(a\omega)                             								\big|^2 d\omega} = a  E_H             							\tag{2.9}
\end{align*}

and
\begin{align*}
  \Delta \Big[ H(a\omega) \Big] = a\Delta H.													 \tag{2.10}
\end{align*}

Thus, the Q-factor of the proposed transform $(2.3)$ is given by

\begin{align*}
  Q = \dfrac{\text{width of the window function}}{\text{centre of the window function}}	= \dfrac{\Delta \big[ H(a\omega) \big] }				{E \big[ H(a\omega) \big]}	= \dfrac{\Delta_H}{E_H}	= \text{constant},				 \tag{2.11}
\end{align*}

which is independent of the unimodular matrix $M = (A,B,C,D:p,q)$ and the scaling parameter $a$. Therefore, the localized time and frequency characteristics of the proposed transform $(2.3)$ are given in the time and frequency windows

\begin{align*}
  \Big[ b+aE_\psi - a\Delta_\psi , b + aE_\psi +a\Delta_\psi \Big] \quad \text{and} \quad \Big[ aE_H - a\Delta_H , aE_H - a				\Delta_H  \Big] ,    																		 \tag{2.12}
\end{align*}

respectively. Hence, the joint resolution of the continuous special affine wavelet transform $(2.3)$ in the time-frequency domain is described by a flexible window having a total spread $4\Delta_\psi\Delta_H$ and is given by

\begin{align*}
  \Big[ b + aE_\psi - a\Delta_\psi , b + aE_\psi + a\Delta_\psi \Big] \times \Big[ aE_H - a\Delta_H , aE_H - a\Delta_H  \Big]     																		 						 \tag{2.13}
\end{align*}

\parindent=0mm\vspace{.1in}
\section{Basic Properties of The Convolution Based Special Affine Wavelet Transform}

In this section we study some fundamental properties of the special affine wavelet transform $(2.3)$ which are similar to the conventional wavelet transform . In this direction, we have the following theorem which assembles some of the basic properties of the proposed transform.

\parindent=0mm\vspace{.1in}
{\bf Theorem 3.1.} {\it For any $f,g \in L^2(\mathbb R)$ and $\alpha,\beta,k \in \mathbb R$, $\mu\in\mathbb R^+$, the continuous special affine wavelet transform defined by $(2.3)$ satisfies the following properties:

\begin{enumerate}
 \item Linearity: \quad $ \mathcal W_{\bf M}^{\psi} \big[ \alpha f + \beta g \big](b,a) =  \alpha \, \mathcal W_{\bf M}^{\psi}    				  \big[ f \big] (b,a) + \beta \, \mathcal W_{\bf M}^{\psi} \big[ g \big](b,a) $ . \\
 \item Translation: \quad $ \mathcal W_{\bf M}^{\psi} \big[ f(t-k)\big](b,a) = \exp \left\{ -i\frac{A}{2B} \Big( 							       \frac{2bk}{a} - \big( 1 + \frac{1}{a^2} \big)k^2 \Big) \right\} \\
                   \mathcal W_{\bf M}^{\psi} \left[ \exp \left\{ i \frac{A}{B}t - k \right\} f(t) \right] (b-k,a)$. \\
 \item Scaling: \quad $ \mathcal W_{\bf M}^{\psi} \big[ f(\mu t) \big](b,a) = \dfrac{1}{\mu} \, \mathcal W_{\bf M}^{\psi} \left[ 			 \exp \left\{ -i \frac{A}{2B}t^2 \big( 1 - \frac{1}{\mu^2} \big) \right\} f(t) \right] (\mu b,\mu a) $.
\end{enumerate}
}

\parindent=0mm\vspace{.1in}
{\it Proof.} The proof of the above theorem is quiet simple, hence is omitted here. \quad\fbox

\parindent=0mm\vspace{.1in}

{\bf Theorem 3.2 (Admissibility Condition).} {\it A given function $\psi \in L^2(\mathbb R)$ is said to be admissible if

\begin{align*}
  C_\psi = \int_{\mathbb R^+} \frac{ \Big| \mathrm S_{\bf M} \left[ \exp \left\{ \frac{i}{2B} \Big( 2zp(a-1) - Az^2 \Big)    				\right\} \overline{\psi(-t)} \right] (a\omega) \Big|^2 }{a} \,da < \infty, \quad a.e.					 \tag{3.1}
\end{align*}
}

\parindent=0mm\vspace{.1in}
{\it Proof.} For any square integrable function $f$, we have

\begin{align*}
   &\int_\mathbb R \int_\mathbb {R^+} \Big| \Big\langle  f, \psi_{b,a}^{\bf M} \Big\rangle \Big|^2 \frac{db \, da}{a^2}  \\
   &\quad = \int_\mathbb R \int_\mathbb {R^+} \left| \exp \left\{ -i \frac{A}{2B} \Big( 1-\frac{1}{a^2} \Big) b^2 \right\}   	   		\left[ f(t)*_{\bf M} \frac{1}{\sqrt{a}} \, \overline{\psi \left( \frac{-t}{a} \right)} \, \exp \left\{ -i \frac{A}			   		{2B} t^2 \right\}  \right](b) \right|^2 \frac{db \, da}{a^2} 	\\
   &\quad = \int_\mathbb R \int_\mathbb {R^+} \left| \mathrm S_{\bf M} \left[ f(t)*_{\bf M} \frac{1}{\sqrt{a}} \, \overline{    			\psi \left( \frac{-t}{a} \right)} \, \exp \left\{ -i \frac{A} {2B} t^2 \right\}  \right] (\omega) \right|^2 \frac{d					\omega \, da}{a^2} 	\\
   &\quad = \int_\mathbb R \int_\mathbb {R^+} \Big| \mathrm S_{\bf M} \big[ f \big] (\omega) \Big|^2 \left| \mathrm S_{\bf M} 				\left[ \exp \left\{ \frac{i}{2B} \Big( 2zp(a-1) - Az^2 \Big) \right\} \overline{\psi(-t)} \right] (a\omega) \right|^2 				\frac{d	\omega \, da}{a} 	\\
    &\quad = \int_\mathbb R \Big| \mathrm S_{\bf M} \big[ f \big] (\omega) \Big|^2 \left\{ \int_\mathbb {R^+} \frac{\left| 					\mathrm S_{\bf M} \left[ \exp \left\{ \frac{i}{2B} \Big( 2zp(a-1) - Az^2 \Big) \right\} \overline{\psi(-t)} \right] (a				\omega) \right|^2}{a} \, da \right\} d\omega  . 												\tag{3.2}
\end{align*}

For $f=\psi$, $(3.2)$ boils down to

\begin{align*}
  &\int_\mathbb R \int_\mathbb {R^+} \Big| \Big\langle  f, \psi_{b,a}^{\bf M} \Big\rangle \Big|^2 \frac{db \, da}{a^2}  \\
   &\quad = \int_\mathbb R \Big| \mathrm S_{\bf M} \big[ \psi \big] (\omega) \Big|^2 \left\{ \int_\mathbb {R^+} \frac{\left| 					\mathrm S_{\bf M} \left[ \exp \left\{ \frac{i}{2B} \Big( 2zp(a-1) - Az^2 \Big) \right\} \overline{\psi(-t)} \right] 			   (a\omega) \right|^2}{a} \, da \right\} d\omega .  											\tag{3.3}
\end{align*}

Since $\psi \in L^2(\mathbb R)$, therefore we conclude that the R.H.S of $(3.3)$ is finite provided

\begin{align*}
  C_\psi = \int_\mathbb {R^+} \frac{\left| \mathrm S_{\bf M} \left[ \exp \left\{ \frac{i}{2B} \Big( 2zp(a-1) - Az^2 \Big) \right			   \}\overline{\psi(-t)} \right] (a\omega) \right|^2}{a} \, da  < \infty, \quad a.e, \, \omega\in\mathbb R,
\end{align*}

which completes the proof of the Theorem $3.2$. \quad\fbox

\parindent=0mm\vspace{.1in}

We are now in a position to derive the orthogonality relation for the proposed transform $(2.3)$. As a consequence of this formula, we shall deduce the resolution of identity for the proposed transform $(2.3)$.

\parindent=0mm\vspace{.1in}
{\bf Theorem 3.3 (Moyal's Principle).} {\it Let $ \mathcal W_{\bf M}^{\psi} \big[ f \big](b,a) $ and $ \mathcal W_{\bf M}^{\psi} \big[ g \big](b,a) $ be the special affine wavelet transforms of $f$ and $g$ belonging to $L^2(\mathbb R)$, respectively. Then, we have

\begin{align*}
  \int_\mathbb {R} \int_\mathbb {R^+} \mathcal W_{\bf M}^{\psi} \big[ f \big](b,a) \, \overline{\mathcal W_{\bf M}^{\psi} \big[ 			g \big](b,a)} \, \frac{db \, da}{a^2} = C_\psi \, \Big\langle  f,g \Big\rangle_{L^2(\mathbb R)}  ,  	\tag{3.4}
\end{align*}

where $C_\psi$ is given by $(3.1)$.
}

\parindent=0mm\vspace{.1in}
{\it Proof.} Applying Proposition $2.4$, we have for any pair of square integrable functions

\begin{align*}
  \mathcal W^\psi_{\bf M} \big[ f \big] (b,a) &= \sqrt{a}\, \exp \left\{ -i \frac{A}{2B} \Big( 1 - \frac{1}{a^2} \Big) b^2    	\right\} \int_\mathbb R \exp \left\{ \frac{i}{2B} \Big( 2a\omega(Dp-Bq) - Da^2\omega^2 \Big) \right\}  \\
  &\quad\times \mathrm S_{\bf M} \big[ f \big](\omega) \, \mathrm S_{\bf M} \left[ \exp \left\{ \frac{i}{2B} \Big( 2tp(a-1) - 				At^2 \Big) \right\} \overline{\psi(-t)}\right](a\omega) \, \mathcal K_{\bf M^{-1}}(t,\omega) \, d\omega
\end{align*}

and
\begin{align*}
  \mathcal W^\psi_{\bf M} \big[ g \big] (b,a) &= \sqrt{a}\, \exp \left\{ -i \frac{A}{2B} \Big( 1 - \frac{1}{a^2} \Big) b^2    	\right\} \int_\mathbb R \exp \left\{ \frac{i}{2B} \Big( 2a\eta(Dp-Bq) - Da^2\eta^2 \Big) \right\}  \\
  &\quad\times \mathrm S_{\bf M} \big[ g \big](\eta) \, \mathrm S_{\bf M} \left[ \exp \left\{ \frac{i}{2B} \Big( 2tp(a-1) - 			At^2 \Big) \right\} \overline{\psi(-t)}\right](a\eta) \, \mathcal K_{\bf M^{-1}}(t,\eta) \, d\eta.
\end{align*}

Consequently, we have

\begin{align*}
  &\int_\mathbb {R} \int_\mathbb {R^+} \mathcal W_{\bf M}^{\psi} \big[ f \big](b,a) \, \overline{\mathcal W_{\bf M}^{\psi} \big[ 			g \big](b,a)} \, \frac{db \, da}{a^2}   \\
  &\quad = \int_{\mathbb{R\times R\times R \times R^+}} \exp \left\{ \frac{i}{2B} \Big( 2a(Dp-Bq)(\omega-\eta)  - Da^2 \big( 				   											\omega^2-\eta^2 \big) \Big) \right\}    \\
  &\hspace{2cm} \times  \mathrm S_{\bf M} \big[ f \big](\omega) \, \mathrm S_{\bf M} \left[ \exp \left\{ \frac{i}{2B} \Big( 									2tp(a-1) - At^2 Big) \right\} \overline{\psi(-t)}\right](a\omega)   \\
  &\hspace{2.5cm} \times \overline{ \mathrm S_{\bf M} \big[ g \big](\eta) \, \mathrm S_{\bf M} \left[ \exp \left\{ \frac{i}{2B} 											\Big( 2tp(a-1) - At^2 \Big) \right\} \overline{\psi(-t)}\right](a\eta) } \\
  &\hspace{3cm} \times \mathcal K_{\bf M^{-1}}(t,\omega) \, \overline{\mathcal K_{\bf M^{-1}}(t,\eta)} \, \frac{dt\,d\omega\,d								\eta\,da}{a}		\\
  &\quad = \int_{\mathbb{R\times R\times R^+}} \exp \left\{ \frac{i}{2B} \Big( 2a(Dp-Bq)(\omega-\eta)  - Da^2 \big( 				   											\omega^2-\eta^2 \big) \Big) \right\}    \\
  &\hspace{2cm} \times  \mathrm S_{\bf M} \big[ f \big](\omega) \, \mathrm S_{\bf M} \left[ \exp \left\{ \frac{i}{2B} \Big( 									2tp(a-1) - At^2 Big) \right\} \overline{\psi(-t)}\right](a\omega)   \\
  &\hspace{2.5cm} \times \overline{ \mathrm S_{\bf M} \big[ g \big](\eta) \, \mathrm S_{\bf M} \left[ \exp \left\{ \frac{i}{2B} 								\Big( 2tp(a-1) - At^2 \Big) \right\} \overline{\psi(-t)}\right](a\eta) } \\
  &\hspace{3cm} \times \left\{ \int_\mathbb R \mathcal K_{\bf M^{-1}}(t,\omega) \, \overline{\mathcal K_{\bf M^{-1}}(t,\eta)} \, 		 dt \right\}  \frac{d\omega\,d\eta\,da}{a}		
  \end{align*}
 \begin{align*}
  &\quad = \int_{\mathbb{R\times R\times R^+}} \exp \left\{ \frac{i}{2B} \Big( 2a(Dp-Bq)(\omega-\eta)  - Da^2 \big( 				   												\omega^2-\eta^2 \big) \Big) \right\}    \\
  &\hspace{2cm} \times  \mathrm S_{\bf M} \big[ f \big](\omega) \, \mathrm S_{\bf M} \left[ \exp \left\{ \frac{i}{2B} \Big( 											2tp(a-1) - At^2 Big) \right\} \overline{\psi(-t)}\right](a\omega)   \\
  &\hspace{2.5cm} \times \overline{ \mathrm S_{\bf M} \big[ g \big](\eta) \, \mathrm S_{\bf M} \left[ \exp \left\{ \frac{i}{2B} 								\Big( 2tp(a-1) - At^2 \Big) \right\} \overline{\psi(-t)}\right](a\eta) } \\
  &\hspace{3cm} \times \delta(\omega - \eta) \,  \frac{d\omega\,d\eta\,da}{a}    \\
  &\quad = \int_\mathbb {R\times R^+} \mathrm S_{\bf M} \big[ f \big](\omega)\, \mathrm S_{\bf M} \big[ g \big]								     (\omega) \left| \mathrm S_{\bf M} \left[ \exp \left\{ \frac{i}{2B} \Big( 2tp(a-1) - At^2 Big) \right\} 							   \overline{\psi(-t)}\right](a\omega)  \right|^2 \frac{d\omega\,da}{a}	   \\
  &\quad = \int_\mathbb {R} \mathrm S_{\bf M} \big[ f \big](\omega)\, \overline{\mathrm S_{\bf M} \big[ g \big](\omega)} \left\{  			 \int_\mathbb {R^+} \frac{\left| \mathrm S_{\bf M} \left[ \exp \left\{ \frac{i}{2B} \Big( 2tp(a-1) - At^2 Big) 						\right\} \overline{\psi(-t)}\right](a\omega)  \right|^2}{a} \right\} d\omega    \\
  &\quad = C_\psi \Big\langle \mathrm S_{\bf M} \big[ f \big](\omega), \mathrm S_{\bf M} \big[ g \big](\omega) \Big								\rangle_{L^2(\mathbb R)}
   \end{align*}
 \begin{align*}
    = C_\psi \Big\langle f,g \Big\rangle_{L^2(\mathbb R)}.   \hspace{11cm}
 \end{align*}

This completes the proof of the Theorem $3.3$. \quad\fbox

\parindent=0mm\vspace{.1in}

{\it Remarks:}
\begin{enumerate}
  \item For $f = g$, Theorem $3.3$ yields the energy preserving relation associated with the special affine wavelet transform $				(2.3)$. i.e;

    \begin{align*}
       \int_\mathbb {R}\int_\mathbb {R^+} \left| \mathcal W^\psi_{\bf M} \big[ f \big] (b,a) \right|^2 \frac{db\,da}{a^2} = C_   		 \psi \, \big\|  f \big\|^2_{L^2(\mathbb R)}.                                				\tag{3.5}
    \end{align*}

  \item  For $C_\psi = 1$, the operator $\mathcal W^\psi_{\bf M}$ becomes an isometry from $L^2(\mathbb R)$ to $L^2(\mathbb {R					\times R^+})$.
\end{enumerate}

\parindent=0mm\vspace{.1in}

{\bf Theorem 3.4 (Inversion Formula).} {\it If $\mathcal W^\psi_{\bf M} \big[ f \big] (b,a)$ is the special affine wavelet transform of an arbitrary function $f \in L^2(\mathbb R)$, then f can be reconstructed as

\begin{align*}
  f(t)= \frac{1}{C_\psi}  \int_\mathbb {R}\int_\mathbb {R^+} \mathcal W^\psi_{\bf M} \big[ f \big](b,a) \, \psi_{b,a}^{\bf M}(t) 					\, \frac{db\,da}{a^2}, \quad a.e.                                            \tag{3.6}
\end{align*}
 }

\parindent=0mm\vspace{.1in}
{\it Proof.} By virtue of Theorem $3.3$, we have

\begin{align*}
   \Big\langle f,g  \Big\rangle &= \frac{1}{C_\psi}  \int_\mathbb {R}\int_\mathbb {R^+} \mathcal W^\psi_{\bf M} \big[ f \big]					(b,a) \, \overline{\mathcal W^\psi_{\bf M} \big[ g \big](b,a)} \, \frac{db\,da}{a^2} 	\\
   &= \frac{1}{C_\psi}  \int_\mathbb {R}\int_\mathbb {R^+} \mathcal W^\psi_{\bf M} \big[ f \big](b,a) \left\{ \int_\mathbb R  				\overline{g(t) } \, \psi_{b,a}^{\bf M}(t) \, dt \right\} \frac{db\,da}{a^2} 	\\
   &= \frac{1}{C_\psi}  \int_\mathbb {R\times R\times R^+} \mathcal W^\psi_{\bf M} \big[ f \big](b,a)\, \psi_{b,a}^{\bf M}(t) \, 		\overline{g(t)} \, \frac{dt\,db\,da}{a^2}  	\\
   &= \frac{1}{C_\psi} \left\langle \int_\mathbb {R}\int_\mathbb {R^+} \mathcal W^\psi_{\bf M} \big[ f \big](b,a) \, \psi_{b,a}				^{\bf M}(t) \, \frac{db\,da}{a^2}, g(t)  \right\rangle .
\end{align*}

Since $g$ is chosen arbitrarily from $L^2(\mathbb R)$, therefore it follows that

\begin{align*}
  f(t)= \frac{1}{C_\psi}  \int_\mathbb {R}\int_\mathbb {R^+} \mathcal W^\psi_{\bf M} \big[ f \big](b,a) \, \psi_{b,a}^{\bf M}(t) 					\, \frac{db\,da}{a^2}, \quad a.e.
\end{align*}

This completes the proof of the Theorem $3.4$. \quad\fbox

\parindent=8mm\vspace{.1in}

As a consequence of the next theorem, we shall demonstrate that the range of the special affine wavelet transform is a reproducing kernel Hilbert space.

\parindent=0mm\vspace{.1in}
{\bf Theorem 3.5 (Characterization of Range).} {\it If $f\in L^2(\mathbb {R\times R^+})$, then $f$ is the special affine wavelet transform of a certain square integrable function if and only if

\begin{align*}
  f \big( b^\prime,a^\prime \big)= \frac{1}{C_\psi} \int_{\mathbb R} \int_{\mathbb R^+} f ( b,a ) \, \Big\langle 				\psi_{b,a}^{\bf M} , \psi_{b^\prime,a^\prime}^{\bf M} \Big\rangle \, \frac{db\,da}{a^2},    				\tag{3.7}
\end{align*}

where $\psi$ satisfies $(3.1)$.
}

\parindent=0mm\vspace{.1in}
{\it Proof.} Let $f$ belongs to the range of the proposed transform $\mathcal W_{\bf M}^{\psi}$. Then, there exists a square integrable function $g$, such that $ \mathcal W_{\bf M}^{\psi} \,g = f $. In order to show that $f$ satisfies $(3.7)$, we proceed as

\begin{align*}
  f \big( b^\prime,a^\prime \big) &=  \mathcal W_{\bf M}^{\psi} \big[ g \big]\big( b^\prime,a^\prime \big)  \\
  	&= \int_\mathbb R g(t)\, \overline{\psi_{b^\prime,a^\prime}^{\bf M} (t)} \, dt    \\
  	&= \frac{1}{C_\psi} \int_{\mathbb R} \left\{ \int_{\mathbb R} \int_{\mathbb R^+} \mathcal W_{\bf M}^{\psi} \big[ g \big] 			   (b,a) \, \psi_{b,a}^{\bf M} (t) \, \frac{db\,da}{a^2} \right\} \overline{\psi_{b^\prime,a^\prime}^{\bf M} (t)} \, dt \\
  	&= \frac{1}{C_\psi} \int_{\mathbb R} \int_{\mathbb R^+} \mathcal W_{\bf M}^{\psi} \big[ g \big](b,a) \left\{ \int_{\mathbb 				R} \psi_{b,a}^{\bf M} (t) \, \overline{\psi_{b^\prime,a^\prime}^{\bf M} (t)} \, dt \right\} \frac{db\,da}{a^2}		\\
  	&= \frac{1}{C_\psi} \int_{\mathbb R} \int_{\mathbb R^+} f ( b,a ) \, \Big\langle \psi_{b,a}^{\bf M} , \psi_{b^\prime,a^					\prime}^{\bf M} \Big\rangle \, \frac{db\,da}{a^2},
\end{align*}

which clearly proves our claim. Conversely, suppose that a square integrable function $f$ satisfies $(3.7)$. Then, we show
that there exist a function $g\in L^2(\mathbb R)$ satisfying $ \mathcal W_{\bf M}^{\psi} \,g = f $. The desired function $g$ is constructed as follows

\begin{align*}
  g(t) =  \frac{1}{C_\psi} \int_{\mathbb R} \int_{\mathbb R^+} f ( b,a ) \, \psi_{b,a}^{\bf M} (t) \, \frac{db\,da}{a^2} .																																\tag{3.8}
\end{align*}

It is straightforward to obtain  $ \big\|  g  \big\| \leq \big\|  f  \big\| < \infty $, i.e; $g\in L^2(\mathbb R)$. Moreover, as a consequence of the well known Fubini-theorem, we have

\begin{align*}
   \mathcal W_{\bf M}^{\psi} \big[ g \big]\big( b^\prime,a^\prime \big) &= \int_{\mathbb R} g(t)\,\overline{\psi_{b^\prime,a^			 \prime}^{\bf M} (t)} \, dt		\\
   &= \frac{1}{C_\psi} \int_{\mathbb R} \left\{ \int_{\mathbb R} \int_{\mathbb R^+} f(b,a) \, \psi_{b,a}^{\bf M} (t) \, \frac{db		  \,da}{a^2} \right\} \overline{\psi_{b^\prime,a^\prime}^{\bf M} (t)} \, dt \\
   &= \frac{1}{C_\psi} \int_{\mathbb R} \int_{\mathbb R^+} f(b,a) \left\{ \int_{\mathbb R} \psi_{b,a}^{\bf M} (t) \, 					 \overline{\psi_{b^\prime,a^\prime}^{\bf M} (t)} \, dt \right\} \frac{db\,da}{a^2} 	\\
   &= \frac{1}{C_\psi} \int_{\mathbb R} \int_{\mathbb R^+} f(b,a) \Big\langle \psi_{b,a}^{\bf M} (t), \psi_{b^\prime,a^\prime}			 ^{\bf M} (t) \Big\rangle \frac{db\,da}{a^2} 	\\
   &= f\big( b^\prime,a^\prime \big).
\end{align*}

This completes the proof of the Theorem $3.5$. \quad\fbox

\parindent=0mm\vspace{.1in}

{\bf Corollary 3.5 (Reproducing Kernel Hilbert space).} {\it For any admissible wavelet $\psi\in L^2(\mathbb R)$, the range of the proposed transform $(2.3)$ is a reproducing kernel Hilbert space in $\psi\in L^2(\mathbb {R\times R^+})$ with the kernel given by

\begin{align*}
  \mathcal K_{\bf M}^{\psi} \big( b,a;b^\prime,a^\prime \big) = \Big\langle \psi_{b,a}^{\bf M}, \psi_{b^\prime,a^\prime}^{\bf M}  		\Big\rangle_{L^2(\mathbb R)}.  																	\tag{3.9}
\end{align*}
}

\parindent=0mm\vspace{.1in}
\section{Special Affine Multi-resolution Analysis}

\parindent=0mm\vspace{.1in}
In this section, we introduce the theory of special affine multi-resolution analysis as it sets the ground for the discrete special affine wavelet transform and the construction of orthogonal special affine wavelets.

\parindent=0mm\vspace{.1in}
{\bf Definition 4.1.} A sequence of closed subspaces $\left\{ V_j^{\bf M} \right\}_{j\in\mathbb Z}$ of $L^2(\mathbb R)$ is called a special affine multi-resolution analysis if the following properties hold:

\parindent=0mm\vspace{.1in}
\begin{enumerate}
  \item[(a)] $ V_j^{\bf M} \subset V_{j+1}^{\bf M} $, for all $j\in\mathbb Z$;
  \item[(b)] $ \cup_{j\in\mathbb Z} \, V_j^{\bf M}$ is dense in $ L^2(\mathbb R)$;
  \item[(c)] $\cap_{j\in\mathbb Z} \, V_j^{\bf M} = \{0\}$;
  \item[(d)] $f(t)\in  V_j^{\bf M}$ if and only if $\exp \Big\{ i\frac{3A}{2B} t^2 \Big\} f(2t)\in V_{j+1}^{\bf M} $, for all 						$j\in\mathbb Z$;
  \item[(e)] There exists a function $\phi(t)$ called the scaling function or the father wavelet such that $\left\{ \phi^{\bf M}					_{0,k}(t)= \phi(t-k)\, \exp \left\{ -\frac{i}{2B} \big( At^2 + Dp^2 - Ak^2 \big) \right\}:k\in\mathbb Z \right\}				  $ is an orthonormal basis of subspace $V^{\bf M}_{0}$.
\end{enumerate}

\parindent=0mm\vspace{.1in}
In the above definition, if we assume that the set of functions $\left\{ \phi^{\bf M}_{0,k}(t):k\in\mathbb Z \right\}$ forms a Riesz basis of $V^{\bf M}_{0}$, then $\phi(t)$  generates a generalized special affine multi-resolution analysis $\left\{ V_j^{\bf M} \right\}_{j\in\mathbb Z}$ of $L^2(\mathbb R)$  such that

\begin{align*}
  \phi^{\bf M}_{j,k}(t)= 2^{j/2} \, \phi \big( 2^j t - k \big) \, \exp \left\{ -\frac{i}{2B} \Big( At^2+Dp^2 - Ak^2 \Big) \right\}
\end{align*}

is an orthonormal basis of $ \left\{ V_j^{\bf M} \right\}_{j\in\mathbb Z} $.

\parindent=0mm\vspace{.1in}
{\it Remarks.}

\begin{enumerate}
  \item When $\bold M = (A,B,C,D:0,0)$, the Definition $(4.1)$ boils down to the definition of multi-resolution analysis 						associated with LCT \cite{JW}.
 \item For $\bold M = (\cos\theta,\sin\theta,-\sin\theta,\cos\theta:0,0)$, the Definition $(4.1)$ reduces to the 						definition of multi-resolution analysis associated with fractional wavelet transform \cite{Fr}.
 \item When $\bold M = (0,1,-1,0:0,0)$, one recovers from the Definition $(4.1)$, the definition of the classical multi-						resolution analysis \cite{MM}.
\end{enumerate}

\parindent=0mm\vspace{.1in}
For every $ j \in \mathbb Z $, we define $ W^{\bf M}_j $ to be the orthogonal compliment of $ V^{\bf M}_j $ in $ V^{\bf M}_{j+1} $. Then, we have

\begin{align*}
  V^{\bf M}_{j+1} = V^{\bf M}_j \oplus W^{\bf M}_j \quad \text{and} \quad W^{\bf M}_\ell \perp W^{\bf M}_k \quad \text{if} \quad 		\ell \neq k.
\end{align*}

It follows that for $ j > J $,

\begin{align*}
   V^{\bf M}_{j} = V^{\bf M}_{J} \oplus \bigoplus_{ \ell = 0 }^{ j - J - 1 } W^{\bf M}_{ j - \ell } ,
\end{align*}

where all these subspaces are orthogonal. By virtue of condition $(b)$ in the Definition $4.1$, this implies

\begin{align*}
  L^2(\mathbb R) = \bigoplus_{j\in \mathbb Z} W^{\bf M}_{j},
\end{align*}

a decomposition of $ L^2(\mathbb R) $ into mutually orthogonal subspaces.

\parindent=0mm\vspace{.1in}
{\bf Theorem 4.2.} {\it If $\phi\in L^2(\mathbb R) $, then the collection $\left\{ \phi^{\bf M}_{0,k}(t):k\in\mathbb Z \right\}$ is a Riesz basis of the space $V^{\bf M}_0$ of $ L^2(\mathbb R) $ if and only if there exists positive constants $A_1$ and $A_2$ such that for all $\omega\in [0, 2\pi B]$, we have

\begin{align*}
  A_1 \leq \sum_{k\in\mathbb Z} \left| \widehat{\phi} \left(\frac{\omega-p}{B} + 2k\pi\right) \right|^2 \leq A_2.  \tag{4.1}
\end{align*}
}

\parindent=0mm\vspace{.1in}
{\it Proof.} For any $f(t)\in V^{\bf M}_0 $, we have

\begin{align*}
  f(t) = \sum_{k\in\mathbb Z}  c\big[k \big]\, \phi^{\bf M}_{0,k}(t), 										\tag{4.2}
\end{align*}

where $c\big[k \big] \in \ell^2(\mathbb Z)$.

\parindent=8mm\vspace{.1in}
By implementing SAFT on both sides of $(4.2)$, we obtain

\begin{align*}
   \mathrm S_{\bf M} \big[ f\big](\omega) &= \int_\mathbb R \sum_{k\in\mathbb Z}  c\big[k \big] \, \phi(t-k)\, \exp \left\{ -				\frac{i}{2B}\Big( At^2+Dp^2-Ak^2 \Big) \right\} \mathcal K_{\bf M}(t,\omega)\, dt		\\
   &= \frac{1}{\sqrt{i2\pi B}} \, \sum_{k\in\mathbb Z}  c\big[k \big] \, \exp \left\{ \frac{i}{2B} \Big( Ak^2 -2\omega(Dp-Bq) + 			D\omega^2 \Big) \right\} 		\\
   &\hspace{2.5cm} \times \int_\mathbb R \phi(t-k)\exp\left\{ \frac{i}{B} t(p-\omega) \right\} dt		\\
   &= \frac{1}{\sqrt{i2\pi B}} \, \sum_{k\in\mathbb Z}  c\big[k \big] \, \exp \left\{ \frac{i}{2B} \Big( Ak^2 + 2k(p-\omega) - 				2\omega(Dp-Bq) + D\omega^2 \Big) \right\} 		\\
   &\hspace{3cm} \times \int_\mathbb R \phi(z)\exp\left\{ -i \left( \frac{\omega-p}{B} \right) z \right\} dz 		\\
   &= \exp\left\{ -\frac{i}{2B} Dp^2 \right\} \widetilde{C}_{\bf M}(\omega) \, \widehat{\phi} \left( \frac{\omega-p}{B} \right), 																								     \tag{4.3}
\end{align*}

\parindent=0mm\vspace{.0in}
where $\widetilde{C}_{\bf M}(\omega)$ is the DT-SAFT of $c\big[k\big]$.

\parindent=8mm\vspace{.1in}
Since $\left| \widetilde{C}_{\bf M}(\omega) \right|$ is a $2\pi B-$periodic function \cite{11}, therefore by invoking Parsevals formula for the SAFT, we have

\begin{align*}
  \big\| f \big\|^2_{L^2(\mathbb R)} &= \Big\| \mathrm S_{\bf M} \big[ f\big](\omega) \Big\|^2_{L^2(\mathbb R)} 		\\
  & = \int_\mathbb R \left| \widetilde{C}_{\bf M}(\omega) \right|^2 \left| \widehat{\phi} \left( \frac{\omega-p}{B} \right) 				\right|^2 d\omega			\\
  &= \sum_{k\in\mathbb Z} \int_{2\pi Bk}^{2\pi B(k+1)} \left| \widetilde{C}_{\bf M}(\omega) \right|^2 \left| \widehat{\phi} 				\left( \frac{\omega-p}{B} \right) \right|^2 d\omega			\\
  &= \sum_{k\in\mathbb Z} \int_{0}^{2\pi B} \left| \widetilde{C}_{\bf M}(\omega + 2k\pi B) \right|^2 \left| \widehat{\phi} 				\left( \frac{\omega-p}{B} + 2K\pi \right) \right|^2 d\omega	 		\\
  &=  \int_{0}^{2\pi B} \left| \widetilde{C}_{\bf M}(\omega) \right|^2 \sum_{k\in\mathbb Z}\left| \widehat{\phi} 							\left( \frac{\omega-p}{B} + 2K\pi \right) \right|^2 d\omega	.											\tag{4.4}
\end{align*}

\parindent=0mm\vspace{.0in}
Now,

\begin{align*}
   &\int_{0}^{2\pi B} \left| \widetilde{C}_{\bf M}(\omega) \right|^2 d\omega 		\\
    &\qquad= \frac{1}{2\pi B} \int_{0}^{2\pi B} c\big[ k\big] \, \overline{c\big[ \ell\big]} \, \exp \left\{ \frac{i}{2B} \Big( 			A\big(k^2-\ell^2 \big) - 2\omega(k-\ell)+2p(k-\ell) \Big) \right\} d\omega
\end{align*}

and

\begin{align*}
  \int_{0}^{2\pi B} \exp\left\{ -\frac{i}{2B} 2\omega(k - \ell) \right\} d\omega &= \int_{0}^{2\pi B} \exp\left\{ \frac{i}{2B} 				2\omega( \ell - k) \right\} d\omega 		\\
  &= B \int_{0}^{2\pi} \exp\big\{ i\eta ( \ell - k) \big\} d\eta  	\\
  &= 2\pi B \, \delta_{k,\ell}.
\end{align*}

Therefore, it follows that

\begin{align*}
  \left\| \widetilde{C}_{\bf M}(\omega) \right\|^2_{L^2[0,2\pi B]}= \int_{0}^{2\pi B} \left| \widetilde{C}_{\bf M}(\omega) 					\right|^2 d\omega = \sum_{k\in\mathbb Z} \Big| c\big[k\big] \Big|^2 = \Big\| c\big[k\big] \Big\|^2_{\ell^2 }.
\end{align*}

Since

\begin{align*}
  0< A_1 \leq \sum_{k\in\mathbb Z} \left| \widehat{\phi} \left(\frac{\omega-p}{B} + 2k\pi\right) \right|^2 \leq A_2<\infty
\end{align*}

and
\begin{align*}
   \Big\| c\big[k\big] \Big\|^2_{\ell^2 } = \left\| \exp \left\{ \frac{1}{2B} ak^2 \right\} c\big[k\big] \right\|^2_{\ell^2 },
\end{align*}

we have

\begin{align*}
  A_1 \left\| \widetilde{C}_{\bf M}(\omega) \right\|^2 = A_1 \, \Big\| c\big[k\big] \Big\|^2 \leq \Big\| \mathrm S_{\bf M} \big[ 	 f\big](\omega) \Big\|^2 \leq A_2 \, \Big\| c\big[k\big] \Big\|^2 = A_2 \left\| \widetilde{C}_{\bf M}(\omega) \right\|^2.
\end{align*}

In more explicit form,we can write

\begin{align*}
 A_1 \, \Big\| c\big[k\big] \Big\|^2 \leq \left\| \sum_{k\in\mathbb Z} c\big[k\big]\,\phi^{\bf M}_{0,k}(t) \right\|^2 \leq A_2 				  \,\Big\| c\big[k\big] \Big\|^2.
\end{align*}

This completes the proof of the Theorem $4.2$.\quad\fbox

\parindent=0mm\vspace{.1in}

{\bf Theorem 4.3 (Orthonormalization Process).} {\it If $\phi \in L^2(\mathbb R)$ and if $\left\{ \phi^{\bf M}_{0,k}(t) :k\in\mathbb Z \right\}$ is a Riesz basis, then $\left\{ \tilde{\phi}^{\bf M}_{0,0}(t) :k\in\mathbb Z \right\}$ is an orthonormal basis of $V^{\bf M}_0$ with

\begin{align*}
  \hat{\tilde{\phi}} \left( \frac{\omega-p}{B} \right) =  \frac{\widehat{\phi} \left( \frac{\omega-p}{B} \right)}{\left(   					\displaystyle\sum_{k\in\mathbb Z} \left| \widehat{\phi} \left( \frac{\omega-p}{B} + 2k\pi \right)\right|^2  \right)^{1/2}} .
\end{align*}
}

\parindent=0mm\vspace{.1in}
{\it Proof.} Since $\phi (t)$ is a scaling function associated with the special affine multi-resolution analysis $\left\{ V_j^{\bf M} \right\}_{j\in\mathbb Z}$ and  $\left\{ \phi^{\bf M}_{0,k}(t) :k\in\mathbb Z \right\}$ is a Riesz basis of $V^{\bf M}_0$ of $L^2(\mathbb R)$, therefore we have

\begin{align*}
 \tilde{\phi}^{\bf M}_{0,0}(t)  = \sum_{k\in\mathbb Z} c\big[k\big] \, \phi(t-k)\exp \left\{ -\frac{i}{2B} \Big( At^2 + Dp^2 - 	Ak^2 \Big) \right\}.																		 \tag{4.5}
\end{align*}

By implementing SAFT on both sides of $(4.5)$, we obtain

\begin{align*}
  \hat{\tilde{\phi}} \Big( \frac{\omega-p}{B} \Big) &= \sum_{k\in\mathbb Z} c\big[k\big] \exp \left\{ \frac{i}{2B} Ak^2 				\right\} \int_\mathbb R \phi(t-k) \exp\left\{ \frac{i}{B} t(p-\omega) \right\} dt		\\
  &= \sum_{k\in\mathbb Z} c\big[k\big]  \exp \left\{ \frac{i}{2B} \Big(Ak^2 + 2k(p-\omega) \Big) \right\} \int_\mathbb R 				\phi(z)\, \exp\left\{ \frac{i}{B} z(p-\omega) \right\} dz		\\
  &= \widehat{E} \Big( \frac{\omega-p}{B} \Big) \, \widehat{\phi} \Big( \frac{\omega-p}{B} \Big),   			 \tag{4.6}
\end{align*}

where $e\big[k\big] = c\big[k\big]\exp \left\{ \frac{i}{2B} Ak^2 \right\}$ and $ \widehat{E} \big( \frac{\omega-p}{B} \big) $ is DTFT of $e\big[k\big]$.

\parindent=8mm\vspace{.1in}
If $\left\{ \tilde{\phi}^{\bf M}_{0,k}(t) = \tilde{\phi}(t-k)\exp\left\{ -\frac{i}{2B} \Big( At^2 +Dp^2 -Ak^2 \Big) \right\} : k\in\mathbb Z \right\}$ is an orthonormal basis of $V^{\bf M}_0$, then from Theorem 4.2, we have

\begin{align*}
  \sum_{k\in\mathbb Z} \left| \hat{\tilde \phi} \left( \frac{\omega-p}{B} + 2k\pi \right) \right|^2=1 		\tag{4.7}
\end{align*}

By combining $(4.6)$ and $(4.7)$, we obtain

\begin{align*}
  \left| \widehat{E} \left( \frac{\omega-p}{B} \right) \right|^2\, \sum_{k\in\mathbb Z} \left| \widehat{\phi} \left( \frac{\omega-		p}{B}+ 2k\pi \right) \right|^2	=	1 ,																    \tag{4.8}
\end{align*}

which yields
\begin{align*}
  \left| \widehat{E} \left( \frac{\omega-p}{B} \right) \right|^2 = \frac{1}{\displaystyle\sum_{k\in\mathbb Z} \left| 					\widehat{\phi} \left( \frac{\omega-p}{B}+ 2k\pi \right) \right|^2}  .										\tag{4.9}
\end{align*}

From $(4.6)$ and $(4.9)$, we obtain

\begin{align*}
\hat{\tilde{\phi}} \Big( \frac{\omega-p}{B} \Big) = \frac{\widehat{\phi} \left( \frac{\omega-p}{B} \right)}{ \left( 					\displaystyle\sum_{k\in\mathbb Z} \left| \widehat{\phi} \left( \frac{\omega-p}{B}+ 2k\pi \right) \right|^2 \right)^{1/2}}.  																													\tag{4.10}
\end{align*}

This completes the proof of the Theorem 4.3. \quad\fbox

\parindent=0mm\vspace{.1in}

{\bf Theorem 4.4.} {\it If $\phi$ is the scaling function of the special affine multi-resolution analysis, then the following equality holds

\begin{align*}
  \left| \widehat{C}_0 \left( \frac{\omega-p}{2B} \right)  \right|^2 + \left| \widehat{C}_0 \left( \frac{\omega-p}{2B} + \pi					\right)\right|^2 = 2.                                                        \tag{4.11}
\end{align*}
}

\parindent=0mm\vspace{.1in}
{\it Proof.} The scaling function $\phi\in V^{\bf M}_0$ implies that $\phi\in V^{\bf M}_1$ and $\Big\{ \phi^{\bf M}_{1,k}(t): k\in\mathbb Z \Big\}$ is an orthonormal basis of $\phi\in V^{\bf M}_1$. Thus the scaling function $\phi$ has the following representation.

\begin{align*}
  \phi^{\bf M}_{0,0}(t) = \sqrt{2} \, \sum_{k\in\mathbb Z} c\big[k\big] \, \phi(2t-k)\exp \left\{ -\frac{i}{2B} \Big( At^2 +Dp^2 	-Ak^2 \Big) \right\},      																	\tag{4.12}
\end{align*}

where $ c\big[k\big] = \Big\langle \phi^{\bf M}_{0,0},\phi^{\bf M}_{1,k} \Big\rangle $ and $\displaystyle\sum_{k\in\mathbb Z} \Big| c\big[k\big] \Big|^2 < \infty $.

\parindent=8mm\vspace{.1in}
By implementing SAFT on both sides of $(4.12)$, we obtain

\begin{align*}
  \hat{\phi} \left( \frac{\omega-p}{B} \right) &= \frac{1}{\sqrt{2}} \, \sum_{k\in\mathbb Z} c^{\bf M}\big[k\big] \exp \left\{ -		\frac{i}{2B} k(\omega-p)\right\} \hat{\phi} \left( \frac{\omega-p}{2B} \right)		\\
  &= \frac{1}{\sqrt{2}} \, \widehat{C}_0 \left( \frac{\omega-p}{2B} \right) \hat{\phi} \left( \frac{\omega-p}{2B} \right),    																											\tag{4.13}
\end{align*}

\parindent=0mm\vspace{.0in}
where $\widehat{C}_0 \left( \frac{\omega-p}{2B} \right)$ is the DTFT of $c^{\bf M}\big[k\big]$ called as low-pass filter with $c^{\bf M}\big[k\big] = c\big[k\big] \exp\left\{ \frac{i}{2B} Ak^2\right\}$.

\parindent=8mm\vspace{.1in}
By the orthonormality of $\phi$, we have

\begin{align*}
   \sum_{k\in\mathbb Z} \left| \widehat{C}_0 \left( \frac{\omega-p}{2B} + k\pi \right) \right|^2 \left| \widehat{\phi} \left( \frac{\omega-p}{2B}+ k\pi \right) \right|^2 = 2.
\end{align*}

\parindent=0mm\vspace{.0in}
Since $ \widehat{C}_0 (\omega) $ is a $2\pi-$periodic function, therefore splitting $k$ into even and odd parts, we obtain

\begin{align*}
    2 &= \left| \widehat{C}_0 \left( \frac{\omega-p}{2B} \right) \right|^2 \sum_{m\in\mathbb Z} \left| \widehat{\phi} \left( 				 \frac{\omega-p}{2B}+ 2m\pi \right) \right|^2   \\
     &\qquad + \left| \widehat{C}_0 \left( \frac{\omega-p}{2B} + \pi \right)  \right|^2 \sum_{m\in\mathbb Z}  \left| 						\widehat{\phi} \left( \frac{\omega-p}{2B}+ 2m\pi + \pi \right)\right|^2 	\\
     &= \left| \widehat{C}_0 \left( \frac{\omega-p}{2B} \right) \right|^2 + \left| \widehat{C}_0 \left( \frac{\omega-p}{2B} + 					\pi \right)  \right|^2.
\end{align*}

This completes the proof of the Theorem 4.4.\quad\fbox

\parindent=0mm\vspace{.1in}

{\bf Theorem 4.5.} {\it For any two functions $\phi$ and $\psi$ belonging to $L^2(\mathbb R)$, the set of functions $\left\{ \phi^{\bf M}_{0,k}(t): k\in \mathbb Z \right\}$ and $\left\{ \phi^{\bf M}_{0,\ell}(t): \ell\in \mathbb Z \right\}$ are biorthogonal, that is

\begin{align*}
   \Big\langle \phi^{\bf M}_{0,k}, \phi^{\bf M}_{0,\ell} \Big\rangle=0, \quad\text{for all}\,\, k,\ell \in \mathbb Z,
\end{align*}

if and only if

\begin{align*}
 \sum_{k\in \mathbb Z} \widehat{\phi} \left( \frac{\omega-p}{B} + 2k\pi  \right) \overline{\widehat{\psi} \left( \frac{\omega-p}		{B} + 2k\pi \right)} = 0 \quad a.e.													\tag{4.14}
\end{align*}
 }
\parindent=0mm\vspace{.1in}
{\it Proof.} We have

\begin{align*}
   &\Big\langle \phi^{\bf M}_{0,k}, \phi^{\bf M}_{0,\ell} \Big\rangle  \\
   &\quad= \left\langle \mathrm S_{\bf M} \Big[ \phi^{\bf M}_{0,k} \Big](\omega), \mathrm S_{\bf M} \Big[ \phi^{\bf M}_{0,\ell} 				\Big](\omega) \right\rangle 	\\
   &\quad= \exp \left\{i \frac{A}{2B}  \big( k^2-\ell^2 \big) \right\} \int_\mathbb R \exp \left\{ \frac{i}{B} (p-\omega)(k-					 \ell) \right\} \widehat{\phi} \left( \frac{\omega-p}{B} \right) \overline{\widehat{\psi} \left( \frac{\omega-p}{B} 					   \right)} d\omega  \\
   &\quad= \exp \left\{i \frac{A}{2B}  \big( k^2-\ell^2 \big) \right\} \sum_{k\in \mathbb Z} \int_{2\pi Bk}^{2\pi B(k+1)}  \exp 				\left\{ \frac{i}{B} (p-\omega)(k-\ell) \right\}       \\
   &\hspace{3cm} \times \widehat{\phi} \left( \frac{\omega-p}{B} \right) \overline{\widehat{\psi} \left( \frac{\omega-p}{B}						\right)} d\omega  \\
   &\quad= \exp \left\{i \frac{A}{2B}  \big( k^2-\ell^2 \big) \right\} \int_{0}^{2\pi B} \exp \left\{ \frac{i}{B} (p-\omega)(k-					\ell) \right\}
   \end{align*}
 \begin{align*}
   &\hspace{3cm} \times \left[ \sum_{k\in \mathbb Z} \widehat{\phi} \left( \frac{\omega-p}{B} + 2k\pi \right) 									\overline{\widehat{\psi} \left( \frac{\omega-p}{B} +2k\pi \right) } \right]d\omega .
\end{align*}

Thus

\begin{align*}
   \Big\langle \phi^{\bf M}_{0,k}, \phi^{\bf M}_{0,\ell} \Big\rangle=0, \quad\text{for all}\,\, k,\ell \in \mathbb Z,
\end{align*}

if and only if

\begin{align*}
 \sum_{k\in \mathbb Z} \widehat{\phi} \left( \frac{\omega-p}{B} + 2k\pi  \right) \overline{\widehat{\psi} \left( \frac{\omega-p}		{B} + 2k\pi \right)} = 0 \quad a.e.
\end{align*}

This completes the proof of the Theorem 4.5. \quad\fbox

\parindent=8mm\vspace{.1in}

We now return to the main  problem of constructing a mother wavelet $\psi(t)$ from a special affine multi-resolution analysis. Suppose that there is a function $\psi$ such that $\left\{ \psi^{\bf M}_{j,k} (t): k\in \mathbb Z \right\}$ is a basis for a space $W^{\bf M}_{j}$. Then $\psi^{\bf M}_{0,0} (t) \in W^{\bf M}_0 \subseteq V^{\bf M}_1 $ and hence there exists a coefficient sequence $\big\{ d\big[k\big]\big\}_{k\in\mathbb Z}$ such that

\begin{align*}
  \psi^{\bf M}_{0,0} (t) = \sum_{k\in\mathbb Z} d\big[k\big] \phi^{\bf M}_{1,k} (t)    		\tag{4.15}
\end{align*}

\parindent=0mm\vspace{.0in}
Proceeding in the same way as in the derivation of $(4.13)$, we have

\begin{align*}
  \hat{\psi} \left( \frac{\omega-p}{B} \right) = \frac{1}{\sqrt{2}} \, \widehat{C}_1 \left( \frac{\omega-p}{2B} \right) \hat{\phi} \left( \frac{\omega-p}{2B} \right),    											\tag{4.16}
\end{align*}

where
\begin{align*}
  \widehat{C}_1 \left( \frac{\omega-p}{2B} \right) = \sum_{k\in\mathbb Z} d^{\bf M}\big[k\big] \exp\left\{ \frac{i}{2B}k(p-\omega) \right\}.																	\tag{4.17}
\end{align*}

Following the technique implemented in the derivation of $(4.11)$, $\widehat{C}_1$ satisfy the following equality

\begin{align*}
  \left| \widehat{C}_1 \left( \frac{\omega-p}{2B} \right)  \right|^2 + \left| \widehat{C}_1 \left( \frac{\omega-p}{2B} + \pi					\right)\right|^2 = 2.                                                        \tag{4.18}
\end{align*}

Substituting $(4.13)$ and $(4.16)$ in $(4.14)$, we obtain

\begin{align*}
  \sum_{k\in\mathbb Z} \widehat{C}_0 \left( \frac{\omega-p}{2B} + k\pi \right) \overline{ \widehat{C}_1 \left( \frac{\omega-p}{2B}+k\pi \right)}  \left| \widehat{\phi} \left( \frac{\omega-p}{2B} + k\pi \right) \right|^2 = 0.
\end{align*}

Since $ \widehat{C}_0 (\omega) $ and $ \widehat{C}_1 (\omega) $ are $2\pi-$periodic function, therefore splitting $k$ into even and odd parts, we obtain

\begin{align*}
    0 &=  \widehat{C}_0 \left( \frac{\omega-p}{2B} \right) \overline{ \widehat{C}_1 \left( \frac{\omega-p}{2B} \right)} \sum_{m				\in\mathbb Z} \left| \widehat{\phi} \left( \frac{\omega-p}{2B}+ 2m\pi \right) \right|^2   \\
     &\qquad +  \widehat{C}_0 \left( \frac{\omega-p}{2B} + \pi \right) \overline{ \widehat{C}_1 \left( \frac{\omega-p}{2B} + \pi 		\right) } \sum_{m\in\mathbb Z}  \left| \widehat{\phi} \left( \frac{\omega-p}{2B}+ 2m\pi + \pi \right)\right|^2 	\\
     &= \widehat{C}_0 \left( \frac{\omega-p}{2B} \right) \overline{ \widehat{C}_1 \left( \frac{\omega-p}{2B} \right)} + 						\widehat{C}_0 \left( \frac{\omega-p}{2B} + \pi \right) \overline{ \widehat{C}_1 \left( \frac{\omega-p}{2B} + \pi 				\right) }    																\tag{4.19}
\end{align*}

The above information can be put in the determinant form as

\begin{align*}
  \begin{vmatrix}
   \widehat{C}_0 \left( \frac{\omega-p}{2B} \right) & \overline{ \widehat{C}_1 \left( \frac{\omega-p}{2B} + \pi \right) } \\ \\
   -\widehat{C}_0 \left( \frac{\omega-p}{2B} + \pi \right) &  \overline{ \widehat{C}_1 \left( \frac{\omega-p}{2B} \right) }    \end{vmatrix}   =0 .                      											\tag{4.20}
\end{align*}

Relations $(4.13), (4.18) $ and $(4.19)$ shows that if $ \left\{ \psi^{\bf M}_{j,k} (t): k\in \mathbb Z \right\} $ is an orthonormal basis for $W^{\bf M}_j$, then we have

\begin{align*}
  MM^* = 2I, 																		\tag{4.21}
\end{align*}

where $*$ denotes the tranjugate of $M$ and $I$ is the identity matrix with

\begin{align*}
M=
  \begin{pmatrix}
    \widehat{C}_0 \left( \frac{\omega-p}{2B} \right) & \widehat{C}_0 \left( \frac{\omega-p}{2B} + \pi \right)  \\ \\
    \widehat{C}_1 \left( \frac{\omega-p}{2B} \right) & \widehat{C}_1 \left( \frac{\omega-p}{2B} + \pi \right)
  \end{pmatrix} .  																	\tag{4.22}										\end{align*}

This shows that if $ \left\{ \psi^{\bf M}_{j,k} (t): k\in \mathbb Z \right\} $ is an orthonormal basis of $W^{\bf M}_j$, then $\widehat{C}_0$ and $\widehat{C}_1$ are the quadrature mirror filters of the special affine wavelets.

\parindent=8mm\vspace{.1in}
Notice that $(4.20)$ can be interpreted as the linear dependence of vectors \\
$ \left( \widehat{C}_0 \left( \frac{\omega-p}{2B} \right), -\widehat{C}_0 \left( \frac{\omega-p}{2B} + \pi \right) \right)$ and $ \left( \overline{ \widehat{C}_1 \left( \frac{\omega-p}{2B} + \pi \right) } \right), \left( \overline{ \widehat{C}_1 \left( \frac{\omega-p}{2B} \right) } \right) $ and hence there exists a function $\lambda$ such that

\begin{align*}
   \widehat{C}_1 \left( \frac{\omega-p}{2B} \right) = \lambda \left( \frac{\omega-p}{2B} \right) \overline{ \widehat{C}_0 \left( \frac{\omega-p}{2B} + \pi \right)} . 												\tag{4.23}
\end{align*}

\parindent=0mm\vspace{.0in}
Substituting $(4.23)$ into $(4.20)$ yields

\begin{align*}
  \lambda \left( \frac{\omega-p}{2B} \right) + \lambda \left( \frac{\omega-p}{2B} + \pi \right) =0 \quad a.e.		\tag{4.24}	
\end{align*}

Letting $\omega = \omega + 2\pi B$, we obtain

\begin{align*}
  \lambda \left( \frac{\omega-p}{2B} + \pi \right) + \lambda \left( \frac{\omega-p}{2B} + 2\pi\right) =0.        	\tag{4.25}
\end{align*}

Equivalently,

\begin{align*}
  \lambda \left( \frac{\omega-p}{2B} \right) =  \lambda \left( \frac{\omega-p}{2B} + 2\pi\right),                \tag{4.26}
\end{align*}

which shows that $\lambda$ is a $2\pi B-$periodic function. Thus, there exists a $2\pi B-$periodic function $\nu$ defined by

\begin{align*}
   \lambda \left( \frac{\omega-p}{2B} \right) = \exp\left\{ - i\frac{\omega}{2B} \right\} \nu \left( \frac{\omega-p}{B} \right). 																											\tag{4.27}
\end{align*}

From $(4.13)$, we have

\begin{align*}
  \widehat{C}_0 \left( \frac{\omega-p}{2B} + \pi \right) &= \sum_{k\in\mathbb Z} c^{\bf M}\big[k\big] \exp \left\{ -\frac{i}{2B}   				k(\omega-p) -ik\pi \right\}   \\
  & \sum_{k\in\mathbb Z} (-1)^k c^{\bf M}\big[k\big]  \exp \left\{ -\frac{i}{2B}k(\omega-p)\right\}.    		\tag{4.28}
\end{align*}

Hence from $(4.23), (4.27)$ and $(4.28)$, it follows that

\begin{align*}
   \widehat{C}_1 \left( \frac{\omega-p}{2B} \right) = \exp\left\{ - i\frac{\omega}{2B} \right\} \nu \left( \frac{\omega-p}{B} 			\right) \sum_{k\in\mathbb Z} (-1)^k \overline{ c^{\bf M}\big[k\big]} \exp \left\{ \frac{i}{2B}k(\omega-p)\right\} . 																														\tag{4.29}
\end{align*}

On comparing $(4.17)$ and $(4.29)$, we obtain

\begin{align*}
  &\sum_{k\in\mathbb Z} d^{\bf M}\big[k\big] \exp\left\{-\frac{i}{2B}k(\omega-p) \right\}		\\
  &\qquad =  \exp\left\{ - i\frac{\omega}{2B} \right\} \nu \left( \frac{\omega-p}{B} \right) \sum_{k\in\mathbb Z} (-1)^k \overline{ c^{\bf M}\big[k\big]} \exp \left\{ \frac{i}{2B}k(\omega-p)\right\} .							  \tag{4.30}
\end{align*}

In particular, if we set $\nu \left( \frac{\omega-p}{B} \right)=1$, then $(4.30)$ yields

\begin{align*}
  &\sum_{k\in\mathbb Z} d^{\bf M}\big[k\big] \exp\left\{-\frac{i}{2B}k(\omega-p) \right\}		\\
   &\qquad =  \exp\left\{ - i\frac{\omega}{2B} \right\}  \sum_{k\in\mathbb Z} (-1)^k \overline{ c^{\bf M}\big[k\big]} \exp \left     	\{ \frac{i}{2B}k(\omega-p)\right\} .							
\end{align*}

Equivalently,

\begin{align*}
    &d^{\bf M}\big[k\big] \exp\left\{ \frac{i}{2B}kp \right\}  \\
    &\qquad = \int_\mathbb R \exp\left\{ - i\frac{\omega}{2B} \right\} \sum_{n\in	\mathbb Z} (-1)^n \overline{ c^{\bf M}\big[n			\big]} \exp \left\{ \frac{i}{2B}n(\omega-p)\right\} \exp\left\{ \frac{i}{2B}k\omega \right\} d\omega
\end{align*}
\begin{align*}
 \Leftrightarrow \quad &d^{\bf M}\big[k\big] \exp\left\{ \frac{i}{2B}kp \right\}  \\
   &\qquad = \sum_{n\in\mathbb Z} (-1)^n \overline{ c^{\bf M}\big[n\big]} \exp \left\{ -\frac{i}{2B}np \right\}  \int_				\mathbb R \exp\left\{ -\frac{i}{2B} \omega(1-n-k) \right\}  d\omega		\\
   &\qquad = \sum_{n\in\mathbb Z} (-1)^n \overline{ c^{\bf M}\big[n\big]} \exp \left\{ -\frac{i}{2B}np \right\} \delta_{1-n-k}\\
   &\qquad = (-1)^{-m+1} \overline{ c^{\bf M}\big[-m+1\big]} \exp \left\{ -\frac{i}{2B}p(-m+1) \right\}.   			\tag{4.31}
\end{align*}

Thus, we obtain a relationship between the coefficients  $d^{\bf M}\big[k\big]$ and $c^{\bf M}\big[k\big]$ as

\begin{align*}
  d^{\bf M}\big[k\big]  = (-1)^{-k+1} \overline{ c^{\bf M}\big[-k+1\big]} \exp \left\{-i \frac{p}{2B} \right\}.																																			\tag{4.32}
\end{align*}

Equivalently,

\begin{align*}
  d\big[k\big]= (-1)^{-k+1} \overline{ c\big[-k+1\big]} \exp \left\{ \frac{i}{2B} \big( A(1-2k)-p \big) \right\}.   \tag{4.33}
\end{align*}

Hence, the unique special affine wavelets can be constructed in this way.

\parindent=8mm\vspace{.1in}
The following is an example showing the construction of special affine Haar wavelet. If we take $\phi(t)= \chi_{[0,1)}(t)$, then the coefficients $d\big[k\big]$ can be obtained according to $(4.33)$. Figs. 3(a), 3(b) and 3(c) describe the constructed special affine Haar wavelets with different parameters as $(A,B,C,D:p,q)=(2,1,1,1:1,1), (A,B,C,D:p,q)=\left(\cos\frac{\pi}{4},\sin\frac{\pi}{4},-\sin\frac{\pi}{4},\cos\frac{\pi}{4} :0,0 \right)$ and $(A,B,C,D:p,q)=(0,1,-1,0:0,0)$ respectively. For the case $(A,B,C,D:p,q)=(2,1,1,1:1,1)$, $d\big[k\big]$ are computed as

\begin{align*}
  d\big[k\big]=
    \begin{cases}
       -\frac{1}{\sqrt{2}}\,(0.4998+0.0131 \, i), \quad k=0  \\
        \frac{1}{\sqrt{2}}\,(0.4998-0.0131 \, i), \quad\quad k=1  \\
        0, \hspace{4.2cm} \text{otherwise}
    \end{cases} 																						\tag{4.34}
\end{align*}

\parindent=0mm\vspace{.0in}
and the corresponding special affine Haar wavelet is obtained as

\begin{align*}
 \psi^{\bf M}_{0,0}(t)=
    \begin{cases}
        -(0.4998+0.0131 \, i) \, \exp\left\{ -\frac{i}{2}\big( 2t^2 + 1 \big) \right\}, \quad 0\leq t<\frac{1}{2}	\\
        (0.4998-0.0131 \, i) \, \exp\left\{ -\frac{i}{2}\big( 2t^2 - 1 \big) \right\}, \quad \frac{1}{2}\leq t \leq 1.	
    \end{cases}																							\tag{3.35}
\end{align*}

Clearly special affine Haar wavelet is more flexible than classical Haar wavelet by choosing different values for unimoludar parametric matrix $\bf M$. For $(A,B,C,D:p,q)=(0,1,-1,0:0,0)$, the special affine Haar wavelet reduces to the classical Haar wavelet as shown in Fig. 3(c).

\section{Construction of special affine multi-resolution analysis from a  Scaling Function}

\parindent=0mm\vspace{.1in}
The main purpose of this section is to construct a special affine multi-resolution analysis by first choosing an appropriate scaling function $\phi $ and obtaining $V^{\bf M}_0$ by taking the linear span of integer translates of $\phi $. The other spaces $V^{\bf M}_j$ can be generated as the scaled versions of $V^{\bf M}_0$.

\parindent=8mm\vspace{.0in}
We begin, with a function $\phi \in L^2(\mathbb R)$ such that

\begin{align*}
 \phi^{\bf M}_{0,0}(t)= \sum_{k\in\mathbb Z} c\big[ k \big] \phi(2t - k)\exp \left \{ -\frac{i}{2B} \Big( At^2+Dp^2-Ak^2 \Big) \right\},                     																	 \tag{5.1}
\end{align*}

\parindent=0mm\vspace{.0in}
where
\begin{align*}
 \sum_{k\in\mathbb Z} \Big| c\big[ k \big] \Big|^2 <\infty \quad \text{and} \quad  0 < A_1 \leq \sum_{k\in\mathbb Z} \left| \widehat{\phi} \Big( \frac{\omega-p}{B} + 2k\pi \Big) \right|^2 \leq A_2 < \infty ,                        \tag{5.2}
\end{align*}

$A_1$ and $A_2$ are constants.

\parindent=8mm\vspace{.0in}
Then, we define $V^{\bf M}_0$ as the closed span of $\big\{ \phi(t-k)\exp\left\{ -\frac{i}{2B} \big( At^2 + Dp^2 - Ak^2 \big)\right\} : k\in \mathbb Z \big\} $ and $V^{\bf M}_j$ as the span of $ \big\{ \phi^{\bf M}_{j,k} :  k \in \mathbb Z \big\}$. The conditions (5.1) and (5.2) are necessary and sufficient to guarantee that $ \big\{ \phi^{\bf M}_{j,k} \big\}_{k \in \mathbb Z }$ is a Riesz basis in each $V^{\bf M}_j$ and that $V^{\bf M}_j$ satisfy the increasing property $V^{\bf M}_j \subset V^{\bf M}_{j+1}$, for all $j \in\mathbb Z$. Moreover, it follows that $V^{\bf M}_j$ satisfy scaling and translating properties of Definition 4.1 also. Now, in order to verify that the ladder of spaces generated by $\phi$ forms a special affine multi-resolution analysis, it is sufficient to show that the following properties also hold:

$$ \bigcap_{j \in\mathbb Z} V^{\bf M}_j = 0 \quad \text{and} \quad \overline{\bigcup_{j \in\mathbb Z} V^{\bf M}_j } = L^2(\mathbb R). $$

\parindent=0mm\vspace{.0in}
The following two theorems verify this.

\parindent=8mm\vspace{.0in}
For $\phi\in L^2(\mathbb R) $, we put

\begin{align*}
  \phi^{\bf M}_{j,k}(t)= 2^j \phi(2^jt-k)\exp\left\{ -\frac{i}{2B} \Big( At^2 +Dp^2 -Ak^2 \Big)  \right\}, \quad k\in \mathbb Z. 																				
\end{align*}

\parindent=0mm\vspace{.1in}

{\bf Theorem 5.1.} {\it Let $\left\{ V^{\bf M}_j \right\}_{j \in \mathbb Z }$ be the family of subspaces satisfying (5.2) and  conditions (a), (d) and (e) of the Definition 4.1. Then $ \bigcap_{j \in\mathbb Z} V^{\bf M}_j = 0 $.
}

\parindent=0mm\vspace{.1in}
{\it Proof.} Since $ \left\{ \phi^{\bf M}_{0,k}: k\in\mathbb Z \right\}$ forms a Riesz basis of $V^{\bf M}_0$. Hence, there exists $A_1,A_2$ with $ 0<A_1\leq A_2<\infty $ such that

\begin{align*}
 A_1 \, \Big\|f \Big\|^2_2 \leq \sum_{k \in\mathbb Z} \Big| \Big\langle f, \phi^{\bf M}_{0,k} \Big\rangle \Big|^2 \leq A_2 \, \Big\|f \Big\|^2_2, \quad f \in V^{\bf M}_0.
\end{align*}

It follows by condition (d) of Definition 4.1 that for all $f\in V^{\bf M}_j $

\begin{align*}
 0< A_1 \, \Big\|f \Big\|^2_2 \leq \sum_{k \in\mathbb Z} \Big| \Big\langle f, \phi^{\bf M}_{j,k} \Big\rangle \Big|^2 \leq A_2 \, \Big\|f \Big\|^2_2<\infty.  																			\tag{5.3}
\end{align*}

\parindent=8mm\vspace{.0in}
Let $ f\in \cap_{j \in\mathbb Z} V^{\bf M}_j $. Then for $\epsilon >0$, there exists a compactly supported continuous function $g$ such that $\| f-g\|_2<\epsilon$. If $\mathcal P^{\bf M}_j$ denotes the orthogonal projection of $ V^{\bf M}_j $, then we have

\begin{align*}
  \Big\| f - \mathcal P^{\bf M}_j g  \Big\|_2 = \Big\| \mathcal P^{\bf M}_j (f - g)  \Big\|_2 \leq \big\| f -g \big\|_2 < \epsilon .
\end{align*}

\parindent=0mm\vspace{.0in}
Therefore, for all $j\in\mathbb Z$, we have

\begin{align*}
  \big\| f \big\|_2 \leq \epsilon + \Big\| \mathcal P^{\bf M}_j  g  \Big\|_2  .
\end{align*}

By using (5.3), we obtain

\begin{align*}
  \Big\| \mathcal P^{\bf M}_j  g  \Big\|_2  \leq A_1^{-1/2} \left[ \sum_{k\in\mathbb Z} \Big| \Big\langle P^{\bf M}_j  g , \phi^{\bf M}_{j,k} \Big\rangle \Big|^2 \right]^{1/2} = A_1^{-1/2} \left[ \sum_{k\in\mathbb Z} \Big| \Big\langle  g , \phi^{\bf M}_{j,k} \Big\rangle \Big|^2 \right]^{1/2} . 														\tag{5.4}
\end{align*}

Since $g$ has a compact support, we can assume that supp $g \subseteq [-N,N], N>0$. We have

\begin{align*}
  \sum_{k\in\mathbb Z} \Big| \Big\langle  g , \phi^{\bf M}_{j,k} \Big\rangle \Big|^2 & \leq 2^j  \sum_{k\in\mathbb Z} \left[ \int_{|t|\leq N} \Big| g(t) \Big| \Big| \phi \big(2^j t - k \big) \Big| dt \right]^2	\\
  &\leq  2^j \, \big\| g \big\|_{\infty} \sum_{k\in\mathbb Z} N \int_{|t|\leq N}  \Big| \phi \big( 2^j t - k \big) \Big|^2 dt \\
  & = N  \, \big\| g \big\|_{\infty} \int_{|y+k|\leq N2^j} \Big| \phi(y) \Big|^2 dy  \\
  &= N  \, \big\| g \big\|_{\infty} \int_{\cup_{k\in\mathbb Z} [k-N2^j, k+N2^j]} \Big| \phi(y) \Big|^2 dy.     \tag{5.5}
\end{align*}

By Dominated Convergence Theorem, the terms on right hand side of (4.3) tends to zero as $j\rightarrow\infty$. In particular, there exists a $j$ such that $ \sum_{k\in\mathbb Z} \left| \left\langle  g , \phi^{\bf M}_{j,k} \right\rangle \right|^2 \leq \epsilon^2 A_1 $.

\parindent=8mm\vspace{.0in}
Hence by (5.3), we have $ \left\| \mathcal P^{\bf M}_j g  \right\|_2 \leq \epsilon $ and therefore, $\big\| f \big\|_2 < 2\epsilon $. But $\epsilon$ is arbitrary, this implies that $\big\| f \big\|_2 =0$. Therefore $f=0$ a.e and hence $\bigcap_{j\in\mathbb Z} V^{\bf M}_j = {0} $.

\parindent=0mm\vspace{.1in}
{\bf Theorem 5.2.} {\it Let $\left\{ V^{\bf M}_j \right\}_{j \in \mathbb Z }$ be the family of subspaces satisfying (5.2) and conditions (a), (d) and (e) of the the Definition 4.1. Assume that $\widehat{\phi}(\omega)$ is continuous at $\omega=0$. Then the following two conditions are equivalent:

\begin{enumerate}
\item $\widehat{\phi}(0)\neq 0$
\item $\overline{\bigcup_{j\in\mathbb Z} V^{\bf M}_j} = L^2(\mathbb R)$.
\end{enumerate}
}

\parindent=0mm\vspace{.1in}

{\it Proof.} Let $ f\in \left( \bigcup_{j\in\mathbb Z} V^{\bf M}_j \right)^\perp $. Then for $\epsilon>0$, there exists a compactly supported $C^\infty-$function $g$ such that $\big\| f-g \big\|_2 < \epsilon $. Then, we have

\begin{align*}
 \Big\| \mathcal P^{\bf M}_j f \Big\|^2_2 = \Big\langle \mathcal P^{\bf M}_j f,\mathcal P^{\bf M}_j f \Big\rangle =\Big\langle f, \mathcal P^{\bf M}_j f \Big\rangle =0
\end{align*}

and
\begin{align*}
 \Big\| \mathcal P^{\bf M}_j g \Big\|^2_2 = \Big\| \mathcal P^{\bf M}_j (f-g) \Big\|^2_2 \leq \Big\| f-g \Big\|^2_2 <\epsilon.
\end{align*}

\parindent=8mm\vspace{.0in}
Since the collection $\left\{ \phi^{\bf M}_{j,k}:k\in\mathbb Z \right\}$ forms a Riesz basis for $V^{\bf M}_j$, therefore there exists $A_1,A_2$ with $0<A_1\leq A_2<\infty $ such that

\begin{align*}
 A_1 \, \Big\|f \Big\|^2_2 \leq \sum_{k \in\mathbb Z} \Big| \Big\langle f, \phi^{\bf M}_{j,k} \Big\rangle \Big|^2 \leq A_2 \, \Big\|f \Big\|^2_2, \quad f \in V^{\bf M}_j.
\end{align*}

\parindent=0mm\vspace{.0in}
In particular, we have

\begin{align*}
  \sum_{k\in\mathbb Z} \Big| \Big\langle  g , \phi^{\bf M}_{j,k} \Big\rangle \Big|^2 \leq A_2 \, \Big\| \mathcal P^{\bf M}_j g \Big\|^2_2 .                      																\tag{5.6}
\end{align*}

Now,
\begin{align*}
  &\sum_{k\in\mathbb Z} \Big| \Big\langle  g , \phi^{\bf M}_{j,k} \Big\rangle \Big|^2  \\
  & \qquad = \dfrac{2^{-j}}{2\pi B} \sum_{k\in\mathbb Z} \left| \int_\mathbb R \mathrm S_{\bf M} \big[ g \big](\omega) \, \overline{\widehat{\phi}\Big( \frac{\omega-p}{2^jB} \Big)} \right.	\\
  &\hspace{3cm} \times \left. \exp \left\{ -\frac{i}{2B} \Big( Ak^2 - 2\omega(Dp-Bq) + D\omega^2 + k(p-\omega)2^{1-j} \Big) \right\} d\omega \right|^2	\\
  & \qquad = \dfrac{2^{-j}}{2\pi B} \sum_{k\in\mathbb Z} \left| \int_0^{2\pi B 2^j} \sum_{\ell\in\mathbb Z} \mathrm S_{\bf M} \big[ g \big](\omega + 2\pi B\ell 2^j) \, \overline{\widehat{\phi}\Big( \frac{\omega-p}{2^jB} + 2\pi\ell \Big)} \right.	\\
  &\hspace{3cm} \times \left. \exp \left\{ -\frac{i}{2B} \Big( Ak^2 - 2\omega(Dp-Bq) + D\omega^2 + k(p-\omega)2^{1-j} \Big) \right\} d\omega \right|^2	\\
  & \qquad = \int_0^{2\pi B 2^j} \left| \sum_{\ell\in\mathbb Z} \mathrm S_{\bf M} \big[ g \big](\omega + 2\pi B\ell 2^j) \, \overline{\widehat{\phi}\Big( \frac{\omega-p}{2^jB} + 2\pi\ell \Big)} \right.	\\
  &\hspace{3cm} \times \left. \exp \left\{ -\frac{i}{2B} \Big( Ak^2 - 2\omega(Dp-Bq) + D\omega^2 + k(p-\omega)2^{1-j} \Big) \right\} d\omega \right|^2	\\
   & \qquad = \sum_{m\in\mathbb Z} \int_\mathbb R \mathrm S_{\bf M} \big[ g \big](\omega) \, \overline{ \mathrm S_{\bf M} \big[ f \big](\omega + 2\pi Bm 2^j) \, \widehat{\phi}\Big( \frac{\omega-p}{2^jB} \Big)} \, \widehat{\phi}\Big( \frac{\omega-p}{2^jB} + 2\pi m \Big) \, d\omega  \\
   & \qquad = \int_\mathbb R \Big| \mathrm S_{\bf M} \big[ g \big](\omega) \Big|^2 \Big| \widehat{\phi} \Big( \frac{\omega-p}{2^jB}  \Big) \Big|^2 d\omega + R(f),
\end{align*}

where
\begin{align*}
  \left|  R(f) \right|= \sum_{m\in\mathbb Z^*} \left| \int_\mathbb R \mathrm S_{\bf M} \big[ g \big](\omega) \, \overline{ \mathrm S_{\bf M} \big[ f \big](\omega + 2\pi Bm 2^j) \, \widehat{\phi}\Big( \frac{\omega-p}{2^jB} \Big)} \, \widehat{\phi}\Big( \frac{\omega-p}{2^jB} + 2\pi m \Big) \, d\omega \right|.
\end{align*}

Since $g\in C^\infty(\mathbb R)$, therefore we may assume that there exists $N>0$, such that

\begin{align*}
  \Big| \mathrm S_{\bf M} \big[ g \big](\omega) \Big| \leq N \Big( 1+|\omega|^2 \Big)^{-3/2}.
\end{align*}

Therefore,

\begin{align*}
  \Big|  R(f) \Big| &\leq N^2 \, \sup_{\omega\in\mathbb R} \left|  \widehat{\phi}(\omega) \right|^2 \sum_{m\in\mathbb Z^*} \int_\mathbb R \Big( 1 + \big| \omega+2^jm\pi \big|^2 \Big)^{-3/2} \Big( 1 + \big| \omega-2^jm\pi \big|^2 \Big)^{-3/2} d\omega	\\
  & \leq N^2 \, \sup_{\omega\in\mathbb R} \left|  \widehat{\phi}(\omega) \right|^2 \sum_{m\in\mathbb R^*} \Big( 1+2^{2j} m^2 \pi^2 \Big)^{-1/2} \int_\mathbb R \Big( 1+|\zeta|^2 \Big)^{-1} d\zeta \\
  & \leq N^\prime \, 2^{-j}.
\end{align*}

Therefore

\begin{align*}
  \int_\mathbb R \Big| \mathrm S_{\bf M} \big[ g \big](\omega) \Big|^2 \Big| \widehat{\phi} \Big( \frac{\omega-p}{2^jB}  \Big) \Big|^2 d\omega \leq A_2 \, \epsilon^2 + N^\prime \, 2^{-j}											\tag{5.7}
\end{align*}

Since $\widehat{\phi}(\omega)$ is continuous at $\omega=0$ with $\widehat{\phi}(0)\neq 0 $, therefore it follows by Dominated Convergence Theorem that the left hand side of (5.7) converges to $\big| \widehat{\phi}(0) \big|^2 \big\|\mathrm S_{\bf M}[ g] \big\|^2_2 $ as $j\rightarrow\infty$. Thus, we have

\begin{align*}
  \Big\|\mathrm S_{\bf M} \big[ g \big] \Big\|^2_2 \leq A_2 \, \epsilon^2 \, \Big| \widehat{\phi}(0) \Big|^{-2} \implies \Big\|\mathrm S_{\bf M} \big[ g \big] \Big\|_2 \leq A_2^{-1/2} \, \epsilon\, \Big| \widehat{\phi}(0) \Big|^{-1}.
\end{align*}

or equivalently
\begin{align*}
  \Big\|\mathrm g \Big\|_2 \leq A_2^{-1/2} \, \epsilon\, \Big| \widehat{\phi}(0) \Big|^{-1}.
\end{align*}

Combining (5.7) with $ \big|\ f-g \big\| \leq \epsilon $, we obtain

\begin{align*}
   \big\| f \big\| \leq \epsilon + \big\| g \big\| \leq \epsilon \left( 1 + A_2^{1/2} \, \Big| \widehat{\phi}(0) \Big|^{-1} \right).
\end{align*}

But since $\epsilon$ is arbitrary, therefore $\big\| f \big\| =0$ and hence $\overline{\bigcup_{j\in\mathbb Z} V^{\bf M}_j} = L^2(\mathbb R)$.

\parindent=8mm\vspace{.1in}
We now assume that $\overline{\bigcup_{j\in\mathbb Z} V^{\bf M}_j} = L^2(\mathbb R)$. Consider $f$ such that $\mathrm S_{\bf M}\big[f\big] = \chi_{[-1,1]}$ so that $\left\| f \right\|^2_2 = \left\| \mathrm S_{\bf M}\big[f\big] \right\|^2_2 = \frac{1}{\pi} $. By condition (a) of Definition 4.1 and our assumption, we have $\left\| f-\mathcal P^{\bf M}f \right\|_2 \rightarrow 0$ as $j\rightarrow \infty$. This implies $\left\| \mathcal P^{\bf M}f \right\|_2 \rightarrow \left\| f \right\|_2$ as $j\rightarrow \infty$. Hence, we have

\begin{align*}
  \Big\| \mathcal P^{\bf M}f \Big\|^2_2 = \left\| \sum_{k\in\mathbb Z} \Big\langle f, \phi^{\bf M}_{j,k} \Big\rangle \phi^{\bf M}_{j,k} \right\|^2_2 =  \sum_{k\in\mathbb Z} \left| \Big\langle f, \phi^{\bf M}_{j,k} \Big\rangle \right|^2,
\end{align*}

\parindent=0mm\vspace{.1in}
as $\left\{ \phi^{\bf M}_{j,k} :k\in\mathbb Z \right\}$ is an orthonormal basis of $V^{\bf M}_j$. From the Plancherel theorem for the SAFT and the fact that $\mathrm S_{\bf M}\big[f\big] = \chi_{[-1,1]}$, we have

\begin{align*}
    &\Big\| \mathcal P^{\bf M}f \Big\|^2_2 	\\
    &\quad=  \sum_{k\in\mathbb Z} \left| \frac{1}{2\pi} \int_\mathbb R \mathrm S_{\bf M}\big[f\big](\omega) \, \frac{2^{j/2}}{\sqrt{2\pi B}} \, \overline{\widehat{\phi}\Big( \frac{\omega-p}{2^jB} \Big) } \right.	\\
     &\hspace{1cm} \times \left. \exp \left\{ -\frac{i}{2B} \Big( Ak^2 - 2\omega(Dp-Bq) + D\omega^2 + k(p-\omega)2^{1-j} \Big) \right\} d\omega \right|^2		\\
     &\quad= \sum_{k\in\mathbb Z} \left| \frac{1}{2\pi} \int_{-1}^{1} \mathrm S_{\bf M}\big[f\big](2^jB\zeta+p) \, \frac{2^{j/2}}{\sqrt{2\pi B}} \, \overline{\widehat{\phi}(\zeta)} \right.	\\
     &\hspace{1cm} \times \left. \exp \left\{ -\frac{i}{2B} \Big( Ak^2 - 2(2^jB\zeta+p)(Dp-Bq) + D(2^jB\zeta+p)^2 + k(-2^jB\zeta)2^{1-j} \Big) \right\} d\zeta \right|^2.
\end{align*}

By Parvesals identity, we have

\begin{align*}
  \Big\| \mathcal P^{\bf M}f \Big\|^2_2 &= \int_{-1}^{1} \frac{2^j}{2\pi B} \,  \left| \mathrm S_{\bf M}\big[f\big](2^jB\zeta+p) \,  \widehat{\phi}(\zeta) \right|^2 d\zeta		\\
  & = \int_{-(1-p)B^{-1}2^{-j}}^{(1-p)B^{-1}2^{-j}} \frac{1}{2\pi B^2} \,  \left| \mathrm S_{\bf M}\big[f\big](\omega)\,\widehat{\phi}\Big( \frac{\omega-p}{2^jB} \Big) \right|^2 d\omega  \\
  &= \int_\mathbb R \frac{1}{2\pi B^2} \left| \mathrm S_{\bf M}\big[f\big](\omega)\,\widehat{\phi}\Big( \frac{\omega-p}{2^jB} \Big) \right|^2 d\omega .
\end{align*}

Since $\Big\| \mathcal P^{\bf M}f \Big\|_2 \rightarrow \frac{1}{\pi}$ as $j\rightarrow \infty$, we have

\begin{align*}
  \lim_{j\rightarrow \infty} \frac{1}{2\pi} \int_\mathbb R \left| \mathrm S_{\bf M}\big[f\big](\omega)\,\widehat{\phi}\Big( \frac{\omega-p}{2^jB} \Big) \right|^2 d\omega  = \frac{B^2}{\pi}.                                   \tag{5.8}
\end{align*}

It follows from the Dominated Convergence Theorem that the left hand side of (5.8) equals to $\left\| \mathrm S_{\bf M}\big[f\big] \right\|^2_2 \left| \widehat{\phi}\Big( \frac{\omega-p}{2^jB} \Big) \right|^2$. Hence $\left| \widehat{\phi}\Big( \frac{\omega-p}{2^jB} \Big) \right|^2 = B^2 \neq 0$. This completes the proof of the Theorem 5.2.\quad\fbox

\newpage
\begin{figure}
  \centering
  \begin{minipage}[b]{0.4\textwidth}
  {\includegraphics[width=\textwidth]{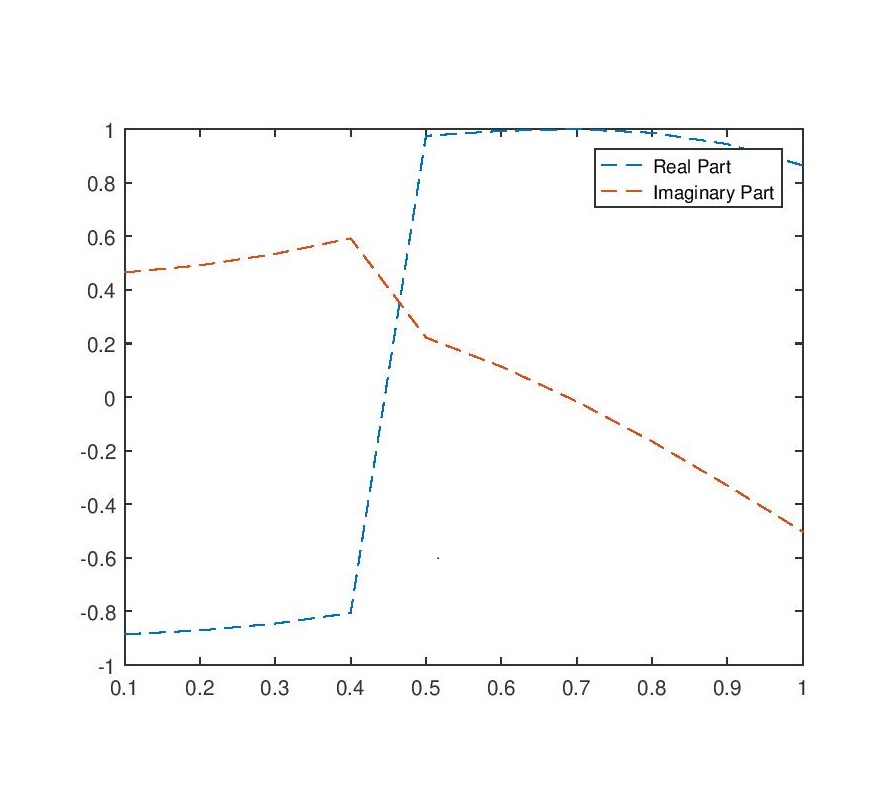}}
  \caption{Special affine Haar wavelet}
  \end{minipage}
  \hfill
\begin{minipage}[b]{0.4\textwidth}
  {\includegraphics[width=\textwidth]{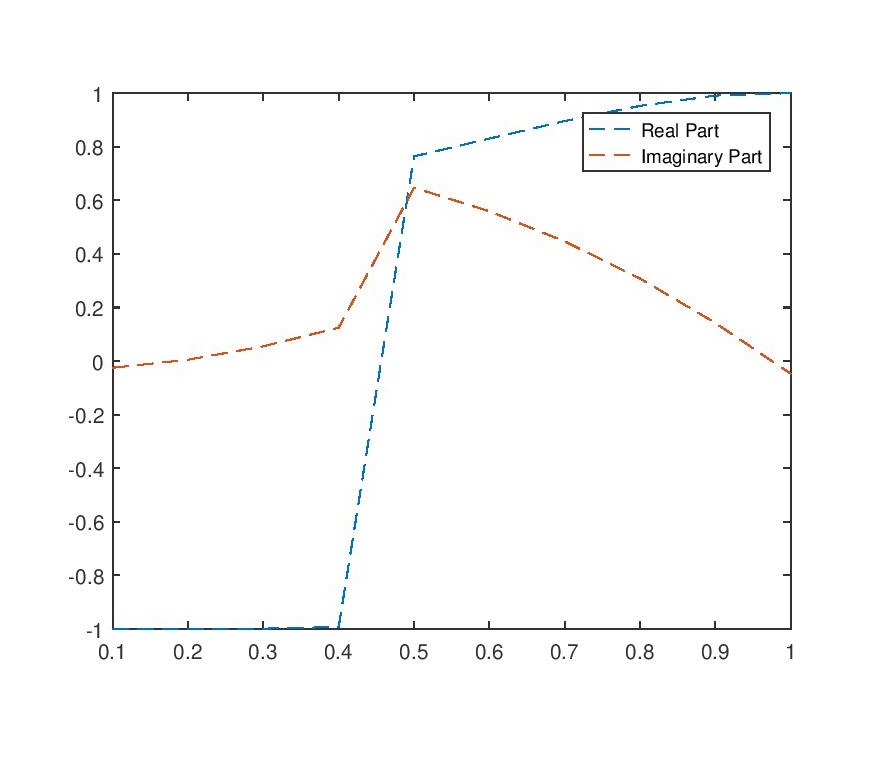}}
  \caption{LCT Haar wavelet}
  \end{minipage}
\end{figure}

\begin{figure}
  \centering
  \begin{minipage}[b]{0.4\textwidth}
  {\includegraphics[width=\textwidth]{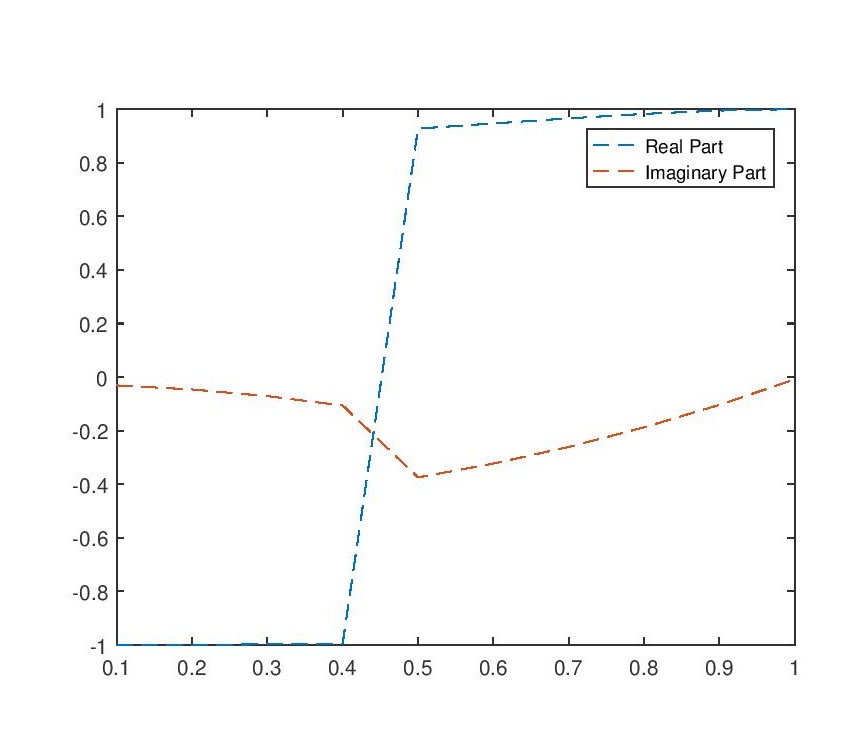}}
  \caption{Fractional Haar wavelet}
  \end{minipage}
  \hfill
\begin{minipage}[b]{0.4\textwidth}
  {\includegraphics[width=\textwidth]{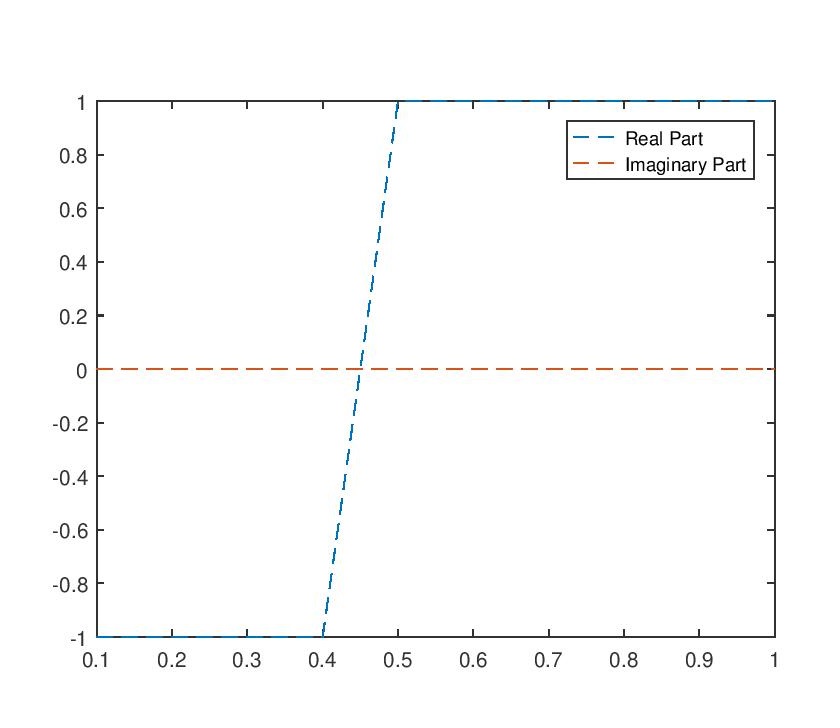}}
  \caption{Classical Haar wavelet}
  \end{minipage}
\end{figure}

\end{document}